\documentclass[onecolumn,3p]{elsarticle}

\usepackage[utf8]{inputenc}

\usepackage{lineno,hyperref}
\modulolinenumbers[50]
 
\journal{Journal of Computational and Applied Mathematics}

\usepackage{cancel}
\usepackage{array}
\usepackage{color}
\usepackage{multirow}
\usepackage{graphicx}
\usepackage{subcaption}
  \usepackage{scrextend}

\usepackage{algorithmicx}
\usepackage{algpseudocode}
\usepackage{algorithm}

\usepackage{bbm}
\usepackage[table]{xcolor}

\usepackage{mathrsfs}
\usepackage{epsfig}
\usepackage{amssymb}
\usepackage{graphicx}
\usepackage{fancyhdr}
\usepackage{float}
\usepackage{epstopdf}

\usepackage{hyperref}

\hypersetup{
     linkcolor = blue,
     citecolor = blue,
     urlcolor = blue,
     anchorcolor = blue
    }

\usepackage{amsthm}

\usepackage{tikz}

\usepackage{pgfplots}
\pgfplotsset{compat=1.11}

\usepackage{xspace}
\usepackage{float}
\usepackage{bbm}
\usepackage{amssymb}
\usepackage{amsmath}
\usepackage{amsthm}
\newcommand{\Sasin}[1]{S^{asin}_#1}

\newcommand{\Pro}[1]{\mathbb{P}\left(#1\right)}

\newcommand{\ie}{\textit{i.e.,}\xspace}
\newcommand{\eg}{\textit{e.g.,}\xspace}

\newtheorem{thm}{Theorem}[section]
\newtheorem{deff}[thm]{Definition}
\newtheorem{lem}[thm]{Lemma}

\newtheorem{exmp}[thm]{Example}

\definecolor{lightgray}{gray}{0.9}

\definecolor{ltgray}{RGB}{200,200,200}
\definecolor{dkgray}{RGB}{150,150,150}

\newcommand{\todo}[1]{}

\begin{document}

\begin{frontmatter}

\title{On testing pseudorandom generators via statistical tests based on the arcsine law\tnoteref{mytitlenote}}
\tnotetext[mytitlenote]{Work supported by NCN Research Grant DEC-2013/10/E/ST1/00359}

\author[uni]{Paweł Lorek}
 \ead{Pawel.Lorek@math.uni.wroc.pl}
\author[gl]{Grzegorz Łoś}
\ead{grzegorz314@gmail.com}
\author[pwr]{Karol Gotfryd}
\ead{Karol.Gotfryd@pwr.edu.pl}
\author[pwr]{Filip Zagórski}
\ead{Filip.Zagorski@pwr.edu.pl}

%
%
%
 \address[uni]{Mathematical Institute, University of Wrocław, pl. Grunwaldzki 2/4, 50-384, Wrocław, Poland}
 \address[gl]{Institute of Computer Science, University of Wrocław, Joliot-Curie 15, 50-383, Wrocław, Poland }
 \address[pwr]{Department of Computer Science, Faculty of Fundamental Problems of Technology, Wrocław University of Science and Technology, Wybrzeże Wyspiańskiego 27, 50-370 Wrocław, Poland}

%
%
%
%
%

\begin{abstract}
Testing the quality of pseudorandom number generators is an important issue. 
Security requirements  become more and more demanding, weaknesses in this matter 
are simply not acceptable. There is a need for an in-depth analysis of statistical tests -- 
one has to be sure that rejecting/accepting a generator as good is not a result of errors in computations or approximations. 
In this paper we propose a second level statistical test based on the arcsine law for random walks.
We provide a Berry-Essen type inequality for approximating the arcsine distribution, what allows us 
to perform a detailed error analysis of the proposed test. 
%
%
%
%
%
\end{abstract}

\begin{keyword}
The arcsine law \sep Random walks \sep Pseudorandom number generator \sep Statistical testing  \sep Second level testing \sep Berry-Esseen type inequality 
\sep Randomness \sep Dyck paths
\end{keyword}

\end{frontmatter}
\linenumbers



 \section{Introduction}\label{sec:intro}
Random numbers are key ingredients in various applications, \eg in cryptography (\eg 
for generating cryptographic keys) or
in simulations (\eg in Monte Carlo methods), just to mention a few. No algorithm can produce truly random numbers.
Instead, pseudorandom number generators (PRNGs) are used. These are deterministic algorithms producing 
numbers which we expect to resemble truly random ones \textsl{in some sense}.
There are two classes of tests used to evaluate PRNGs, theoretical 
and statistical ones. Theoretical tests examine the intrinsic 
structure of a given generator, the sequence does not necessarily need to be 
generated. Two classical examples 
are the \textsl{lattice test}~\cite{MARSAGLIA1972} and the 
\textsl{spectral test} described in~\cite{Knuth2} (Section 3.3.4). 
See also~\cite{LEcuyer92} for a description of some standard tests from 
this class.  This category of tests is very specific to each family of generators
\eg some are designed only for linear congruential generators. 
On the other hand, the second class of tests -- empirical tests --  are
conducted on a sequence generated by a PRNG and require no knowledge of how it 
was produced. The main goal of these tests is to check if the sequence of numbers 
$\mathbf{U} = (U_1, U_2, \ldots, U_n)$ (or bits, depending on the actual 
implementation) produced by a PRNG
has properties similar to those of a sequence 
generated truly at random.  These tests try to find statistical evidence 
against the null hypothesis $\mathcal{H}_0$ stating that the 
sequence is a sample from independent random variables with uniform distribution.
Any function of a finite number of uniformly distributed random variables, whose 
(sometimes approximate) distribution under hypothesis $\mathcal{H}_0$ is known, can
be used as a statistical test.
Due to the  popularity and significance of the problem, a variety of testing procedures have been developed in recent years.
Such statistical tests aim at detecting 
various deviations in generated sequences, what allows for
revealing flawed PRNGs  producing predictable output.
Some of the procedures encompass classical tools from statistics like the Kolmogorov-Smirnov test or 
the Pearson's chi-squared test, which are used for comparing the theoretical and empirical distributions of appropriate
statistics calculated for a PRNG's output. 
It is also possible to adapt tests of normality like the Anderson-Darling or Shapiro-Wilk tests for appropriately transformed
pseudorandom sequences.
These methods exploit the properties of
sequences of i.i.d.
random variables. Based on the original sequence $\mathbf{U}$ returned by
the examined PRNG we are able to obtain
realizations of  random variables with known theoretical distributions. 
Some examples of probabilistic laws
used in practice in this kind of tests can be found \eg in \cite{Knuth2}. They include such procedures like 
the gap test, the permutation test and the coupon collector's test, just to name a few (see \cite{Knuth2} for a more detailed treatment).
These methods have also the advantage that they implicitly test 
the independence of the generator's output.
The main issue with such methods is that a single statistical test looks only at some specific property that holds for
sequences of truly random numbers.  Hence, for practical purposes bundles of diverse tests are created.
Such a test bundle  consists of a series of individual procedures based on various stochastic laws from probability theory.
A PRNG is then considered as good if the pseudorandom sequences it 
produces pass all tests in a given bundle. Note that formally 
it proves nothing,  but it increases the confidence in the simulation results.
Thus, they are actually tests for \textsl{non-randomness}, as pointed out in 
\cite{PareschiRS07}.
Some examples of such test suites are Marsaglia's Diehard Battery of Tests of Randomness from 1995, 
Dieharder developed by Brown et al. (see \cite{Dieharder}), \texttt{TestU01} implemented by L'Ecuyer and Simard
(see \cite{TestU01,TestU01_guide}) and NIST Test Suite \cite{NISTtests}.
The last one, designed by the National Institute of Standard and Technology, is currently
considered as one of the state of the art test bundles. It is often used 
for the preparation of many formal certifications or approvals.

A result of a single statistical test is typically given in the form of a 
$p$-value, which, informally speaking, represents the probability that a perfect 
PRNG would produce ``less random'' sequence than the sequence being tested w.r.t. the 
used statistic. 
We then reject $\mathcal{H}_0$ if $p<\alpha$, where $\alpha$ is the \textsl{significance level} 
(usually $\alpha = 0.01$) and accept $\mathcal{H}_0$ if $p\geq \alpha$.
Such an approach is usually called \textsl{one level} or \textsl{first level} test.
Although the interpretation of a single $p$-value has a clear statistical explanation,
it is not quite obvious how to interpret the results of a test bundle, \ie of multiple tests.
Under $\mathcal{H}_0$ the distribution of $p$-values is uniform. However,
in a test bundle several different tests are applied to the same output of a PRNG,
hence the results are usually correlated.
The documentation of the NIST Test Suite includes some clues on 
how to interpret the results of their bundle 
(Section 4.2  in \cite{Rukhin2010}),
but in the introduction it is frankly stated: 
``It is up to the tester to determine the correct interpretation of the test results''. 

%




To disclose flaws of PRNGs, a very long sequence is often required.
In such situations, the applicability of a statistical test can be limited (depending on the test statistic) by the memory size of the computer.
An alternative  approach is to use a so-called \textsl{two level} (a term used \eg in \cite{LEcuyer92}) 
 or \textsl{second level} 
(a term used \eg in \cite{PareschiRS07,Pareschi2007}) test. In this approach we take into account 
several results from the same test over disjoint sequences generated by a PRNG. 
We  obtain several $p$-values which  are uniformly distributed under $\mathcal{H}_0$,
what is tested by \eg some goodness-of-fit test (with potentially different 
level of significance -- NIST   suggests to use $0.0001$ -- obtaining new ``final'' $p$-value).
The authors in \cite{Ecuyer2002_Sparse} observed that 
this method  may be comparable to a first level test  in terms of the power of a test 
(informally speaking, it represents the probability of observing ``less random'' sequence than 
the sequence being tested under an \textsl{alternative} hypothesis $\mathcal{H}_1$, see  \cite{Ecuyer2002_Sparse}
for details), but often it produces much more 
%
%
%
\textsl{accurate} results, as shown in \cite{PareschiRS07}.
Roughly speaking, the accuracy is related to the ability,
given a non-random PRNG, of recognizing its sequences as non-random 
(for details see \cite{PareschiRS07}). We will follow this approach.

In the second level approach one has to take under consideration the approximation errors in the
computation of a $p$-value.
For example, in a first level test one usually calculates a $p$-value of a statistic which -- under $\mathcal{H}_0$ -- 
is \textsl{approximately} normally distributed. The approximation comes then from  the central limit theorem, which 
lets us substitute the distribution of a given sum with the standard normal distribution. 
These errors in calculations of individual $p$-values may accumulate, resulting
in an error of a $p$-value in
a second level test,
thus making the test \textsl{not reliable}.
Following \cite{PareschiRS07} we say that the second level test is \textsl{not reliable} when, due to errors or approximations 
in the computation of $p$-values (in the first level), the distribution of $p$-values is not uniform under $\mathcal{H}_0$.
Fortunately, this approximation error can be bounded using the Berry-Essen inequality and the final error 
of a second level test 
can be controlled (see  \cite{PareschiRS07,Pareschi2007} for a detailed example based on the binary matrix rank test).
The influence of
approximations on the computation of $p$-values in a second level test was also considered in \cite{Matsumoto2002,Leopardi2009}.
In this article we present a statistical test based on the arcsine law, in which at some point we approximate a distribution
of some random variable with the arcsine distribution. We provide a Berry-Essen type 
inequality which upper bounds the approximation error, what allows us to control the  reliability of our second level test.

 \medskip\par

An interesting approach for testing PRNGs was presented by Kim et al. in \cite{Kim2008}. The concept of their tests is based on the
properties of a random walk (the gambler's ruin algorithm) on the cyclic group $\mathbb{Z}_n=\{0,\ldots,n-1\}$ with 0  being an 
absorbing state -- more 
 precisely, on the time till absorption.
The authors in \cite{Kim2008} propose three different variants of the test. The general idea of the basic procedure 
is the following. For some fixed $p \in (0,1)$ and $x \in \mathbb{Z}_n$, the output $\mathbf{U} = (U_i)$ of a PRNG is treated as 
numbers from the unit interval and used to define a random walk starting in $x$ such that if $U_i < p$ and the process is in state $s$,
then it moves to $(s+1) \mod n$, otherwise it moves to $(s-1) \mod n$. The aim of this test is to 
compare the theoretical and the empirical 
distributions of the time to absorption in 0 when starting at $x$. Based on the values of testing statistic, the PRNG is then either 
accepted or rejected.  The authors reported some ``hidden defects'' in the widely used Mersenne Twister 
generator. However, one has to be very careful when dealing with randomness. It seems like re-seeding a PRNG with a fixed seed is 
an error which can lead to wrong conclusions. The  criticism was raised by  Ekkehard and Grønvik in \cite{Ekke10},
where the authors also showed that the \textsl{properly} performed  tests of Kim et al. \cite{Kim2008} do not reveal any defects in
the Mersenne Twister PRNG.
Recently, the authors in \cite{Lorek2017a} have proposed another gambler's ruin based procedure for testing PRNGs.
In their method they exploited formulas for winning probabilities for arbitrary  sequences $p(i)$ and $q(i)$, (\ie the winning and losing probabilities 
depend on the current fortune) which are the parameters  of the algorithm.

In recent years a novel kind of testing techniques has been introduced for more
careful verification of generators.
The core idea of this class of methods is based on an observation that the binary sequence $(B_i)$ produced by a PRNG,
after being properly rescaled, can be interpreted as an one-dimensional random walk $(S_n)_{n \in \mathbb{N}}$ with 
$S_k = \sum_{i=1}^{k} X_i$, where $X_i = 2 B_i - 1$. For random walks defined by truly random binary sequences a wide range of 
statistics have been considered over the years and a variety of corresponding stochastic laws have been derived
(see \eg \cite{Feller1}). For a good PRNG we may expect that its output will behave like $S_n$. Hence, the following idea comes to
mind: choose some probabilistic  law that holds for truly random bit sequences and compare the theoretical distribution of the 
corresponding statistic with the empirical distribution calculated for $m$ sequences produced by 
a given  PRNG in $m$ 
independent experiments. This comparison can be done \eg by computing the $p$-value of an appropriate test statistic under the 
null hypothesis that the sequence generated by this PRNG is truly random. 

Another concept named \emph{statistical distance based
testing} was suggested in \cite{Wang2015}. It relies on calculation of some statistical distances like
\eg total variation  distance  between the  theoretical and empirical distributions for 
considered characteristics and rejecting a PRNG  if the distance exceeds some threshold.
We will also follow this approach, indicating the corresponding threshold.
In \cite{Wang2015} the authors derive their test statistics from the law of iterated logarithm 
for random walks (the procedure is called the LIL test). The proposed by us 
procedure uses similar methodology and is based on the arcsine law.  
We made the code publicly available, see \cite{ArcsineTest_github}.
It includes the arcsine law based as well as the law 
of iterated logarithm based statistical tests,
the implementation of many PRNGs (more than described in this article) 
including the \textsf{Flawed} generator (see Section \ref{subsec:experiments_flawed}) and the seeds 
we used.

\paragraph{Organization of the paper}
In the following Section \ref{sec:rw} we  define a general notion of a PRNG 
and recall the aforementioned stochastic laws  for random walks. 
The testing method along with the error analysis is  described in Section \ref{sec:arcsine_tester}.
The concise report on experimental results (including the \textsf{Flawed} generator 
introduced in Section \ref{sec:flawed}) is given in 
Section \ref{sec:experiments}. In Section \ref{sec:takashima} we mention other 
implementations of the tests based on the arcsine law.  We conclude in Section \ref{sec:concl}.


 \section{Pseudorandom generators and stochastic laws for  random walks}\label{sec:rw}
 \subsection{Pseudorandom generators}
The intuition behind \textsl{pseudorandom number generator} is clear. However, let us give  a strict definition roughly following Asmussen and Glynn \cite{Asmussen2007}.
\begin{deff}
 \label{def:PRNG}
 A \textsl{Pseudorandom number generator} (PRNG) is a 5-tuple $<E,V,s_0,f,g>$, where $E$ is a finite state space, $V$ is a set of values,
  $s_0\in E$ is a so-called
 \textsl{seed}, \ie an initial state in the  sequence $(s_i)_{i=0}^\infty$, a function $f:E\to E$ describes the transition between consecutive states $s_n=f(s_{n-1})$ 
 and $g: E\to V$ maps the generator's state into the output. 
\end{deff}
Usually $V=(0,1)$ or $V=\{0,1,\ldots,M-1\}$ for some $M\in\mathbb{N}$, the latter one is used throughout the paper.
%
Recall that LCG (linear congruential generator) is a generator which updates its state according 
to the formula $s_n=(a s_{n-1}+c) \bmod M$. Thus, it is defined by three integers: a modulus $M$, a multiplier $a$,
and an additive constant $c$. 
In the case $c=0$, the generator is called MCG (multiplicative congruential generator).
For a detailed description of some commonly used PRNGs see the surveys \cite{Kroese11,LEcuyer2017,Niederreiter_QuasiMonte}
or the book \cite{Woyczynski98}.
\smallskip\par  
It is clear that both the input and the output of a random number generator can be viewed as a finite sequence of bits.
For a PRNG to be considered as good, the output sequences should have some particular property, namely each
returned bit has to be generated independently with equal probability of being 0 and 1.
We say that the sequence of bits is \textsl{truly random} if it is a realization of a  Bernoulli process with success probability $p = \frac{1}{2}$. 

\smallskip\par 
Given a PRNG $G$ returning integers from the set $V$, we may obtain a pseudorandom
binary sequence with any given length using the following simple procedure. Namely, as long as the bit sequence $s$
is not sufficiently long, generate the next pseudorandom number $a$ and append its binary representation
(on $\left\lceil{\log_2 M}\right\rceil$ bits) to the current content of $s$. In the ideal model with $G$ being truly random number
generator, such algorithm produces truly random bit sequences provided that $M$ is a power of 2. Indeed, for $M = 2^{k}$
there is one to one correspondence between $k$-bit sequences and the set $V$. Hence, if each number is generated
independently with uniform distribution on $V$, then each combination of $k$ bits is equally likely and therefore
each bit of the output sequence is independent and equal to 0 or 1 with probability $\frac{1}{2}$.

However, this is not true for $M \neq 2^{k}$. It is easy to observe that in such a case the generator is more 
likely to output 0s and the generated bits are no longer independent. 
Thus, rather than simply outputting the bits of $a$, one may instead take $d$ first bits from the binary representation of $\frac{a}{M}$ for some fixed $d$.
Such a method has the advantage that it can be easily adopted for an underlying generator
returning numbers from the unit interval, what is common for many PRNG implementations.

\subsection{Stochastic laws for random walks}
 Let $(B_i)_{i\geq 0}$ be a \textsl{Bernoulli process} with a parameter $p\in(0,1)$, \ie a sequence of independent random variables with identical distribution $P(B_1=1)=1-P(B_1=0)=p$.
A good PRNG should \textsl{behave} like a generator of Bernoulli process with $p=1/2$ (what we assume from now on). It will be, however, more convenient to consider the following
transformed process
\begin{equation}
\label{eq:xi}
	X_i=2B_i-1, \quad S_0=0, \quad S_k=\sum_{i=1}^k X_i, \ k=1,\ldots
\end{equation}
The sequence $X_i$ is $\{-1,+1\}$-valued, the process $(S_n)_{n\in\mathbb{N}}$ is called a random walk.

\paragraph{The law of iterated logarithm}
Of course $|S_n|\leq n$. However, large values of $|S_n|$ occur with small probability and the values of $S_n$ are in practice in a much narrower range than $[-n,n]$.
\textsl{The weak} and \textsl{the strong law of large numbers} imply that ${S_n\over n}\stackrel{P}{\to}0, \quad \textrm{and even } {S_n\over n}\stackrel{a.s.}{\to}0,$
where $\stackrel{P}{\to}$ denotes the \textsl{convergence in probability} and $\stackrel{a.s.}{\to}$ denotes the \textsl{almost sure convergence}.
Thus, the deviations of $S_n$ from 0 grow much slower than linearly. On the other hand \textsl{the central limit theorem} states that 
${S_n\over n}\stackrel{D}{\to} \mathcal{N}(0,1)$ (where $\stackrel{D}{\to}$ denotes the convergence in distribution),
what is in some sense a lower bound on fluctuations of $S_n$ --  they will leave the interval $[-\sqrt{n},\sqrt{n}]$ since we have 
$\limsup_{n\to\infty} {S_n\over\sqrt{n}}=\infty$
(implied by 0-1 Kolmogorov's Law, see \eg Theorem 5.1 in \cite{Gut2005}). It turns out that the fluctuations can be estimated more exactly.

\begin{thm}[The law of iterated logarithm, \cite{Khintchine1924}, cf. also Chapter VIII.5 in \cite{Feller1}]
 \label{thm:lil}
 For a random walk $S_n$ we have 
 $$\begin{array}{llll}
  & \displaystyle \Pro{\liminf_{n\to\infty} {S_n\over \sqrt{2n\log\log n}}=-1} &  = 1,\\ 
 & \displaystyle  \Pro{\limsup_{n\to\infty} {S_n\over \sqrt{2n\log\log n}}=+1} & =1.
 \end{array}
$$ 
\end{thm}
Thus, to normalize $S_n$ dividing by $n$ is \textsl{too strong} and dividing by $\sqrt{n}$ is \textsl{too weak}. The fluctuations of $S_n$ from 0 grow 
proportionally to $\sqrt{2n\log\log n}$. 

\begin{figure}[H]
\centering
\includegraphics[scale=0.4]{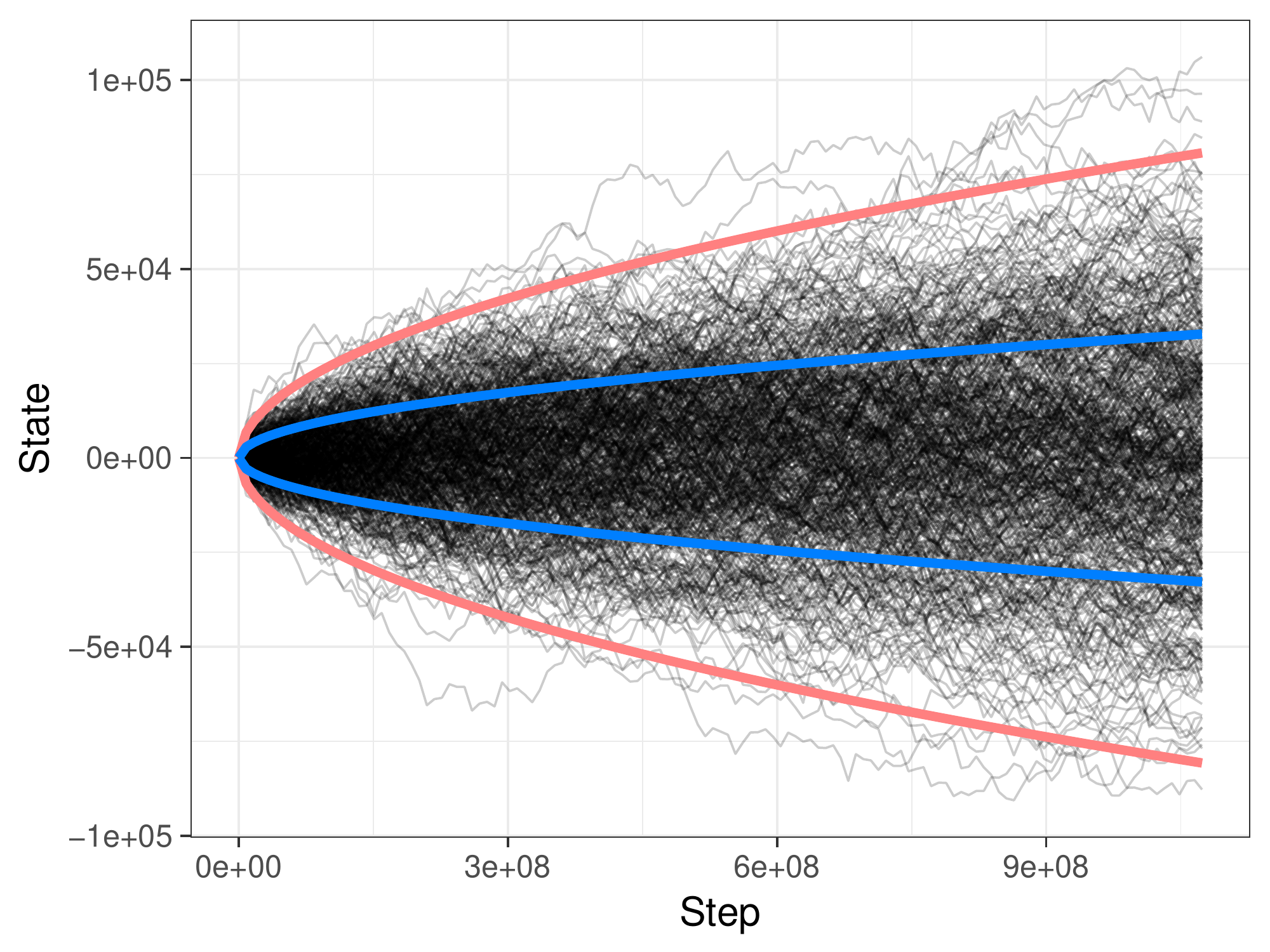}
\caption{500 trajectories of random walks of length $2^{30}$. Blue plot: $\pm \sqrt{n}$, red plot: $\pm \sqrt{2n\log\log n}$ }\label{fig:500iter}
\end{figure}

To depict the law of iterated logarithm, we 
took $500$ output   sequences $\mathcal{B}^1,\ldots,\mathcal{B}^{500}$ from the  Mersenne Twister MT19937 generator,
each initialized with a random seed taken from \url{http://www.random.org},
where each output $\mathcal{B}^j=(B^j_1,\ldots,B^j_n), B^j_i\in\{0,1\}, j=1,\ldots,500, i=1,\ldots,n$ was of length $n=2^{30}$.
In Figure \ref{fig:500iter} we presented these 500 trajectories $(k, S^j_k), j=1,\ldots,500, k=0,\ldots,n$, where $S^j_k=\sum_{i=1}^k (2B_i^j-1)$.
Each trajectory is depicted by a single polyline. The darker the image the higher the density of trajectories. We can see that $\pm \sqrt{2n\log \log n}$ 
roughly corresponds to the fluctuations of $S_n$. However, few trajectories after around billion steps are 
still outside $[-\sqrt{2n\log\log n},\sqrt{2n\log\log n}]$.  
The law of iterated logarithm tells us that for \textsl{appropriately} large $n$ the trajectories 
will not leave $[-\sqrt{2n\log\log n},\sqrt{2n\log\log n}]$ with probability $1$,   what is not the case 
in Figure \ref{fig:500iter}. It means that $n$ must be much larger than $2^{30}$.   
\par 
One could think that the following is a good test for randomness: fix some number, say $100$, and classify the considered PRNG as good
if the difference between the number of ones and zeros never exceeds $100$. The large difference may suggest that 
zeros and ones have  different probabilities of occurrence. However, the law of iterated logarithm tells us that this reasoning is wrong.
Indeed, we \textsl{should} expect some fluctuations and the absence of them means that a PRNG does not produce bits which can be considered random.
This property of random walks was used by the authors in \cite{Wang2015} for designing a novel method of testing
random number generators.

There is yet another interesting property. Define $S_n^{lil}={S_n\over \sqrt{2n\log\log n}}$. 
The law of iterated logarithm  implies that $S_n^{lil}$ does not converge pointwise to any constant. However, it converges to 0
in probability. Let us fix some small $\varepsilon>0$. For almost all $n$, with an arbitrary high 
probability $p<1$ the process $S_n^{lil}$ will not leave $(-\varepsilon,\varepsilon)$.
On the other hand, this tells us that the process will be outside this interval infinitely many times.
This apparent contradiction shows how can our intuition be unreliable on phenomena taking place at infinity.

\paragraph{The arcsine law}
The observations described   previously  imply that \textsl{averaging every} $S_n$, it will spend half of its time above the $x$-axis and half of its time below.
However, the typical situation is counter-intuitive (at first glance): typically the random walk will either spend most of its time above or most of its time below the $x$-axis.
This is expressed in the Theorem \ref{twr:arcsine} below (for reference see \eg \cite{Feller1}).
Before we formulate the theorem, let us first introduce some notations. 
For a sequence  $X_1, X_2, \ldots$, as defined in (\ref{eq:xi}), let
\begin{equation}\label{eq:dk}
	D_k = \mathbbm{1}\left(S_k > 0 \vee S_{k-1}>0 \right), k=1,2,\ldots,
\end{equation}
where  $\mathbbm{1}(\cdot)$ is the indicator function.
$D_k$ is equal to 1 if the number of ones exceeds the number of zeros either at step $k$ or at step $k-1$, and 0 otherwise (in a case of ties, \ie
$S_k = 0$, we look at the previous step letting $D_k = D_{k-1}$).
In other words, $D_k = 1$ corresponds to the situation in which the line segment of the trajectory of the random walk
between steps $k-1$ and $k$ is above the $x$-axis. 
\begin{thm}[The arcsine law]\label{twr:arcsine}

	Let $(B_i)_{i\geq 0}$ be a Bernoulli process.    Define $X_i=2B_i-1$ and $L_{n}=\sum_{k=1}^{n} D_k$ ($D_k$ is given in (\ref{eq:dk})).
	 For $x \in (0,1)$ we have
  \begin{equation*}
		\begin{split}
			\Pro{L_{n}\leq x \cdot n} & \xrightarrow[n \to \infty]{} \frac{1}{\pi} \int_0^x \frac{dt}{\sqrt{t(1-t)}}  =\frac{2}{\pi} \arcsin \sqrt{x}~.\\
		\end{split}
	\end{equation*}
\end{thm}

The probability $\Pro{L_{n}\leq x \cdot n}$ is the chance that the random walk was above the $x$-axis for at most $x$ fraction of the time.
The limiting distribution is called the \textsl{arcsine distribution}. Its density function is given by $f^{asin}(t)={\frac{1}{\pi} \sqrt{t(1-t)}}$
and the cumulative distribution function (cdf) is $F^{asin}(t) = \frac{2}{\pi} \arcsin \sqrt{t}$. 
The shape of the pdf $f^{asin}(t)$ clearly indicates that the fractions of time spent above and below $0$-axis are more likely
to be unequal than close to each other.
%


 \section{Testing PRNGs based on the arcsine law}\label{sec:arcsine_tester}

In this Section we will show how to exploit the theoretical properties of random walks from the preceding discussion
to design a practical routine for testing PRNGs. We describe our approach based on 
the arcsine law which we employ for experimental evaluation of several commonly used generators (the results are presented in Section \ref{sec:experiments}). 
We also perform an error analysis of the proposed testing procedure, providing corresponding bounds on the approximation errors. Finally, we make
some remarks on the reliability of our second level test.

\subsection{The arcsine law based testing}
\label{subsec:arcsineTests}

The general idea of tests is the following. Take a sequence of bits 
generated by PRNG, rescale them as in  (\ref{eq:xi})
and compare the empirical distribution of 
\begin{equation*}
 \Sasin{n} = \frac{1}{n} \sum_{k=1}^n D_k \in[0,1]
\end{equation*}
(a fraction of time instants at which ones prevail zeros)
with its theoretical distribution assuming that truly random numbers were generated.
In terms of hypothesis testing: 
given the null hypothesis $\mathcal{H}_0$ that the bits in the sequence were generated independently and uniformly at random
(vs. $\mathcal{H}_A$: that the sequence was not randomly 
generated), 
the distribution of $\Sasin{n}$
%
 follows the arcsine law (Theorem \ref{twr:arcsine}),
\ie we can conclude that  \textsl{for large} $n$ we have 
\begin{equation}\label{eq:prob_sasin}
\Pro{\Sasin{n} \leq x | \mathcal{H}_0} \approx \frac{1}{\pi} \int_0^x \frac{dt}{\sqrt{t(1-t)}} =  \frac{2}{\pi}\arcsin(\sqrt{x})
\end{equation}
(we will be more specific on ``$\approx$'' in Section \ref{sec:err_analysis}).
We follow the \textsl{second level} testing approach (cf. \cite{PareschiRS07,Pareschi2007}), \ie we take into account several results from the same test over 
different sequences. To test a PRNG we generate $m$ sequences of length $n$ each,
thus obtaining  $m$ realizations of the variable $\Sasin{n}$. Denoting by $\Sasin{{n,j}}$ the value of $j$-th 
simulation's result (we call them a \textsl{basic tests}), we then calculate the
corresponding $p$-values
$$p_j = \Pro{\Sasin{{n}}>\Sasin{{n,j}} | \mathcal{H}_0}= 1-\frac{2}{\pi}\arcsin\left(\sqrt{\Sasin{{n,j}}}\right), j=1,\ldots,m$$
Under $\mathcal{H}_0$ the distribution of $p_j, j=1,\ldots,m$, should be uniform on $[0,1]$.
We fix some partition of  $[0,1]$ and count the number of $p$-values within each interval.
In our tests we will use an $(s+1)$-element partition 
$\mathcal{P}_s = \{ P_1,   \ldots, P_{s+1}\}$, 
where 
\begin{equation*}
\begin{split}
  P_1 &= \left[0, \frac{1}{2s}\right),\\
  P_i &= \left[\frac{2i-3}{2s}, \frac{2i-1}{2s} \right),\ \ \ 2 \leq i \leq s,\\
  P_{s+1} &= \left[1- \frac{1}{2s},1\right].
\end{split}
\end{equation*}
Now we define the measures $\mu_m$ (the uniform measure on $\mathcal{P}_s$), $\nu_n $
(the empirical  measure on $\mathcal{P}_s$),
 $E_i$  (the expected number of $p$-values within $P_i$) and 
 $O_i$ (the number of observed $p$-values within $P_i$).
 For $1 \leq i \leq s+1$ let
\begin{align*}
\mu_m \left( P_i \right) & =  
\left\{
\begin{array}{llll}
 {1\over s} & \mathrm{if } \ i\in \{2,\ldots, s\} \\[12pt]
 {1\over 2s} & \mathrm{if } \ i\in \{1,s+1\} 
\end{array}\right. 
, \quad E_i=
\left\{
\begin{array}{llll}
 {m\over s} & \mathrm{if } \ i\in \{2,\ldots, s\} \\[12pt]
 {m\over 2s} & \mathrm{if } \ i\in \{1,s+1\} 
\end{array}\right. 
\\[12pt]
\nu_m \left( P_i \right) &= \frac{|\{ j:\ p_j \in  P_i, 1 \leq j \leq m\}|}{m},
 \quad O_i=m \cdot \nu_m(P_i). 
\end{align*}

%
\noindent
We perform the Pearson's goodness-of-fit test, which uses the following 
test statistic
\begin{equation*}
\label{eq:chi_asin2}
T ^{asin} = \sum_{i=1}^{s+1} \frac{(O_i - E_i)^2}{E_i}
 = m \cdot\sum_{\mathcal{A}\in \mathcal{P}_s} {(\mu_m(\mathcal{A}) -\nu_m(\mathcal{A}))^2\over 
 \mu_m(\mathcal{A})}.
\end{equation*}
Under the null hypothesis, $T^{asin}$ 
has approximately $\chi^2(s)$ distribution.  
We calculate the corresponding $p$-value 
$$p_{\chi^2} = \Pro{X>T^{asin}},$$
where $X$ has a $\chi^2(s)$ distribution. 
Large values of $T^{asin}$ -- and thus small values of $p_{\chi^2}$ -- let  us suspect 
that a given  PRNG is not good. 
Typically, we reject $\mathcal{H}_0$ (\ie we consider the 
test \textsl{failed}) if $p_{\chi^s}<\alpha,$ 
where $\alpha$ is a predefined \textsl{level of significance}
(for a second level test we use $\alpha=0.0001$, as suggested by NIST).
Note that the probability of rejecting 
$\mathcal{H}_0$  when the sequence is generated by a perfect random 
generator (so-called \textsl{Type I error}) is exactly $\alpha$.
\par 
Another approach relies on the \textsl{statistical distance 
based testing}, which is the technique presented in \cite{Wang2015}.
We consider the statistic 
\begin{equation*}
\label{eq:dtv}
 d_{tv}^{asin} =  {1\over 2} \sum_{\mathcal{A} \in \mathcal{P}_s} 
 |\mu_m(\mathcal{A}) -\nu_m(\mathcal{A})|\in[0,1],
\end{equation*}
\ie a total variation distance between the theoretical 
distribution $\mu_m$ and the empirical 
distribution $\nu_m $. Similarly, large values of 
$d_{tv}^{asin}$ indicate that a given  PRNG is not good. 
Concerning \textsl{Type I error} we will make use 
of the following lemma (see Lemma 3 in \cite{devroye1983}
or its reformulation, Lemma 1 in \cite{Berend2012}).
\begin{lem}\label{lem:devroy}
 Assume $\mathcal{H}_0$ and consider the partition $\mathcal{P}_s$. Then,  for all $\varepsilon \geq \sqrt{20 (s+1)/m} $ we have 
 $$\Pro{2d_{tv}^{asin} > \varepsilon | \mathcal{H}_0} \leq 3 exp\left( -{m\varepsilon^2\over 25}\right).$$
\end{lem}

\medskip\par 

To summarize, for a given PRNG we generate $m$ sequences of length $n$ each. 
and we choose $s$ (and thus the partition $\mathcal{P}_s$). We then
calculate $d^{asin}_{tv}$
 and $T^{asin}$ together with its $p_{\chi^2}$-value. 
We specify the thresholds for $p_{\chi^2}$-value   and 
$d^{asin}_{tv}$ indicating whether the test \textsl{failed} or not 
(the details are presented in Section \ref{sec:experiments}). We denote 
the described procedure as the ASIN test.

\medskip\par\noindent
\textbf{Remark}. Note that the described procedure for calculating 
$T^{asin}$ and $d_{tv}^{asin}$ is equivalent to the following one.
Instead of calculating $p$-values of $\Sasin{{n,j}}$,
we could directly count the number of $\Sasin{{n,j}}$ falling into each interval $P_i, i=1,\ldots,s+1$ 
and compare the empirical distribution with the theoretical one.
To be more precise, for    $1 \leq i \leq s+1$ let 
\begin{align*}
\mu'_m \left( P_i \right) & =  \Pro{\Sasin{n} \in P_i}, \quad E_i=m \cdot \mu'_m(P_i),\\ 
\nu'_m \left( P_i \right) &= \frac{|\{ j:\ \Sasin{{n,j}} \in  P_i, 1 \leq j \leq m\}|}{m},
 \quad O_i=m \cdot \nu'_m(P_i),
\end{align*}
where $\Pro{\Sasin{n} \in P_i}=F^{asin}(b)-F^{asin}(a)$ for $P_i=[a,b]$.
Then statistics $T^{asin}$ and $d_{tv}^{asin}$ can be rewritten as
$$
 T^{asin} =   m \cdot\sum_{\mathcal{A}\in \mathcal{P}_s} {(\mu'_m(\mathcal{A}) -\nu'_m(\mathcal{A}))^2\over 
 \mu'_m(\mathcal{A})}, \qquad  
 d_{tv}^{asin} =  {1\over 2} \sum_{\mathcal{A}\in \mathcal{P}_s} 
 |\mu'_m(\mathcal{A}) -\nu'_m(\mathcal{A})|.$$
This technique was presented in \cite{Wang2015} (for the total variation and few other 
distances) and this is how our implementation of the ASIN test \cite{ArcsineTest_github} calculates 
the statistics.

\medskip\par 
 
We could also calculate just one $p$-value of the statistic $S_{n'}^{asin}$ 
for a longer sequence (say, for $n'=n\cdot m$) -- \ie perform a  \textsl{first level}
test. However, as mentioned in Section \ref{sec:intro},
the \textsl{second level} approach produces 
more accurate results (roughly speaking, the accuracy is related to the ability,
given a non-random PRNG, of recognizing its sequences as non-random, 
see details in \cite{PareschiRS07}). 

%


\medskip\par 

It is worth noting  that the following approach can be
applied when $d^{asin}_{tv}$ or $p_{\chi^2}$ are 
\textsl{slightly} outside the acceptance region (\eg if $p_{\chi^2}\in(10^{-4}, 
10^{-2}$), what suggests rejecting $\mathcal{H}_0$, but is not a \textsl{strong evidence}).
Namely, double the length of the sequence, take a new output from the PRNG and apply the test again.
Repeat the procedure (at most some predefined number of times) until the evidence is strong enough (\eg $p_{\chi^2}<10^{-4}$)
or $\mathcal{H}_0$ is accepted (\eg $p_{\chi^2}>0.01$). This method, 
called \textsl{``automation of statistical tests on randomness''}, was proposed 
and analyzed in \cite{Haramoto2009}.

\medskip\par 
\subsection{Error analysis}\label{sec:err_analysis}
\subsubsection{Bounding errors in approximating $p$-values in basic tests}\label{sec:err_analysis_pval}

In this subsection we will show a bound on the approximation error  in (\ref{eq:prob_sasin}).
Recall that $F^{asin}(x)={2\over \pi} \arcsin\sqrt{t}$.

 \begin{lem}\label{lem:main_error}
 Fix a partition $\mathcal{P}_{s}, s\geq 2$ and an even $n\geq 2$. Let $F_n$ be the cdf of the empirical distribution of $\Sasin{n}$ under $\mathcal{H}_0$
  (stating that the bits $B_1,\ldots,B_n$ were generated uniformly at random), \ie 
  $F_n(x)=\Pro{\Sasin{n} \leq x | \mathcal{H}_0}$. Then we have 
  
%
 \begin{equation*}
 \sup_{x\in[0,1]}|F_n(x)- F^{asin}(x)|\leq {C\over n}, \quad C= {4\over 3\pi} \left(2 -  \frac{3}{2s}\right)   \left(\frac{4s^2}{2s-1}\right)^{\frac{3}{2}}.  
\end{equation*}

%

\end{lem}
\noindent\textsl{Proof}.
We will show that for fixed $a$ and $b$ such that $0\leq a<b\leq 1$ we have 
$$
\left| \Pro{\Sasin{n} \in (a,b)} - \frac{1}{\pi} \int_a^b \frac{dt}{\sqrt{t(1-t)}} \right| \leq {C\over n}. 
$$
Let us assume that   $n=2\mathfrak{n}$.
Let $p_{2k,2\mathfrak{n}}$ denote the probability that during  $2k$ steps in the first $2\mathfrak{n}$ steps the 
random walk was above 0-axis,
\ie
$p_{2k,2\mathfrak{n}}=\Pro{L_{2\mathfrak{n}}=2k}$.
The classical results on a simple random walk state that 
\begin{equation}\label{eq:p2k2n}
	p_{2k,2\mathfrak{n}}={2k\choose k}{2(\mathfrak{n}-k) \choose \mathfrak{n}-k} 2^{-2\mathfrak{n}}.
\end{equation}
The standard proof of Theorem \ref{twr:arcsine} (see, \eg  Chapter XII.8 in \cite{Feller2}) shows that $p_{2k,2\mathfrak{n}}$ converges to  
$d_{k,\mathfrak{n}} = {1\over \pi\sqrt{k(\mathfrak{n}-k)}}$. In the following, we will bound the difference
 $|p_{2k,2\mathfrak{n}}-d_{k,\mathfrak{n}}|$.
We will use a version of Stirling's formula 
stating that for each $\mathfrak{n}$ there exists $\theta_\mathfrak{n}$, $0<\theta_\mathfrak{n}\leq 1$, such that 
\begin{equation}\label{eq:n!}
	\mathfrak{n}!=\sqrt{2\pi \mathfrak{n}} \left({\mathfrak{n}\over e}\right)^\mathfrak{n} \exp \left({\theta_\mathfrak{n}\over 12\mathfrak{n}}\right).
\end{equation}
Plugging (\ref{eq:n!}) into each factorial appearing in (\ref{eq:p2k2n}) we have\par
 \begin{equation*}
p_{2k,2\mathfrak{n}}= \frac{1}{\pi \sqrt{k(\mathfrak{n}-k)}} \exp\left( \frac{\theta_{2k} - 4\theta_k}{24k} + \frac{\theta_{2(\mathfrak{n}-k)} - 4\theta_{\mathfrak{n}-k}}{24(\mathfrak{n}-k)} \right).
\end{equation*}
Thus, we get 
\[ \frac{p_{2k,2\mathfrak{n}}}{d_{k,\mathfrak{n}}} \leq \exp\left( \frac{1}{24k} + \frac{1}{24(\mathfrak{n}-k)} \right) =
\exp\left( \frac{\mathfrak{n}}{24k(\mathfrak{n}-k)} \right) \]
and
\[ \frac{p_{2k,2\mathfrak{n}}}{d_{k,\mathfrak{n}}} \geq \exp\left( \frac{-4}{24k} + \frac{-4}{24(\mathfrak{n}-k)} \right) = 
\exp\left( -\frac{\mathfrak{n}}{6k(\mathfrak{n}-k)} \right). \]
For any $x$ it holds that $1-e^{-x}\leq x$ and for $x\in[0,1.25]$
we have that $e^x-1\leq 2x$. Note that $\frac{\mathfrak{n}}{24(\mathfrak{n}-k)}\leq 1.25$, what is 
equivalent to $\mathfrak{n}\geq {30k^2\over 30k-1}$, what holds for any $0<k<\mathfrak{n}$.
Hence, 
\begin{equation*}
\begin{array}{llllll} 
  p_{2k,2\mathfrak{n}} - d_{k,\mathfrak{n}} & \leq  & \displaystyle d_{k,\mathfrak{n}} \left(\exp\left( \frac{\mathfrak{n}}{24k(\mathfrak{n}-k)} \right) - 1\right) & \leq & \displaystyle  d_{k,\mathfrak{n}}\frac{\mathfrak{n}}{12k(\mathfrak{n}-k)}, \\[10pt]
  d_{k,\mathfrak{n}} - p_{2k,2\mathfrak{n}} & \leq &  \displaystyle d_{k,\mathfrak{n}} \left( 1 -\exp\left( -\frac{\mathfrak{n}}{6k(\mathfrak{n}-k)} \right) \right) & \leq &  \displaystyle  d_{k,\mathfrak{n}}\frac{\mathfrak{n}}{6k(\mathfrak{n}-k)},\\
 \end{array}
\end{equation*}
what implies
\[ |p_{2k,2\mathfrak{n}} - d_{k,\mathfrak{n}}| \leq  d_{k,\mathfrak{n}}\frac{\mathfrak{n}}{6k(\mathfrak{n}-k)} = 
\frac{\mathfrak{n}}{6\pi\left( k(\mathfrak{n}-k) \right)^{\frac{3}{2}}}. \]
Fix  $\delta > 0$ and assume furthermore that $\delta \leq \frac{k}{\mathfrak{n}} \leq 1 - \delta$. 
The function  $k \mapsto \left( k(\mathfrak{n}-k) \right)^{3/2}$ achieves the 
minimum value at the endpoints of the considered interval, thus 
\[ |p_{2k,2\mathfrak{n}} - d_{k,\mathfrak{n}}| \leq \frac{\mathfrak{n}}{6\pi\left( \delta \mathfrak{n} 
(\mathfrak{n}- \delta \mathfrak{n}) \right)^{\frac{3}{2}}} = \frac{1}{6\pi \mathfrak{n}^2 
\left( \delta(1- \delta) \right)^{\frac{3}{2}}}.  \]
We will estimate the approximation error in (\ref{eq:prob_sasin}) in two steps. 
First, take two numbers $a, b$ such that $\delta \leq a < b \leq 1 - \delta$.
We have
 
\begin{align}
	\label{eq:approxError1}
		& \left| \sum_{a \leq \frac{k}{\mathfrak{n}} \leq b} p_{2k,2\mathfrak{n}} -  \sum_{a \leq \frac{k}{\mathfrak{n}} \leq b}  d_{k,\mathfrak{n}} \right| \leq 
				\sum_{a \leq \frac{k}{\mathfrak{n}} \leq b} \left|  p_{2k,2n} -  d_{k,\mathfrak{n}} 
				\right|\leq \sum_{a \leq \frac{k}{\mathfrak{n}} 
				\leq b} \frac{1}{6\pi \mathfrak{n}^2 \left( \delta(1- \delta) \right)^{\frac{3}{2}}} 
				\nonumber \\[8pt]
		& = \frac{\lceil{b\mathfrak{n} - a\mathfrak{n}}\rceil}{6\pi \mathfrak{n}^2 
		\left( \delta(1- \delta) \right)^{\frac{3}{2}}}
		 \leq \frac{b - a}{3\pi \mathfrak{n} \left( \delta(1- \delta) \right)^{\frac{3}{2}}} 
		 \leq \frac{1}{3\pi \mathfrak{n} \left( \delta(1- \delta) \right)^{\frac{3}{2}}}=:\eta. \nonumber
\end{align}
%
%
%
The second kind of errors in probability estimates given by (\ref{eq:prob_sasin}) is caused by approximating the sum by an integral.
Let us consider an arbitrary function $f$ differentiable in the interval $(a,b)$. Split $(a,b)$ into subintervals
of length $\frac{1}{\mathfrak{n}}$ and let $x_k$ be an arbitrary point in the interval containing $\frac{k}{\mathfrak{n}}$.
Denote by $M_k$
and $m_k$ the maximum and the minimum value of $f$ on that interval, respectively. Using the Lagrange's mean value theorem we obtain
\begin{align*}
   & \left| \int_a^b f(x) dx - \sum\limits_{a \leq \frac{k}{\mathfrak{n}} \leq b}\frac{1}{\mathfrak{n}}f(x_k) \right|
				\leq \sum\limits_{a \leq \frac{k}{\mathfrak{n}} \leq b}\frac{1}{\mathfrak{n}}(M_i - m_i)   = \sum\limits_{a \leq \frac{k}{\mathfrak{n}} \leq b}\frac{1}{\mathfrak{n}^2}|f'(\xi_i)| 
				\leq \sum\limits_{a \leq \frac{k}{\mathfrak{n}} \leq b}\frac{1}{\mathfrak{n}^2} \sup_{a \leq x \leq b} |f'(x)| \\[8pt]
   & = \frac{\lceil{b\mathfrak{n}-a\mathfrak{n}}\rceil}{\mathfrak{n}^2} \sup_{a \leq x \leq b} |f'(x)| \leq \frac{2(b-a)}{\mathfrak{n}} \sup_{a \leq x \leq b} |f'(x)|. 
\end{align*}
%
For $f(x) = \frac{1}{\pi \sqrt{x(1-x)}}$ we have $f'(x) = \frac{2x-1}{2 \pi (x(1-x))^{3/2}}$ and
$\frac{1}{\mathfrak{n}}f\left(\frac{k}{\mathfrak{n}} \right) = d_{k,\mathfrak{\mathfrak{n}}}.$ Hence, in the considered interval 
$(a,b) \subseteq (\delta, 1-\delta)$ we have
\begin{align*}
  \left| \int_a^b f(x) dx - \sum\limits_{a \leq \frac{k}{\mathfrak{n}} \leq b}d_{k,\mathfrak{n}} \right|  
	        & \leq \frac{2}{\mathfrak{n}} \sup_{\delta < x < 1-\delta} |f'(x)|  
            = \frac{1 - 2\delta}{\pi \mathfrak{n} (\delta(1-\delta))^{\frac{3}{2}}}=:\kappa.
\end{align*}
We also have 

\begin{align*}
\left| \int_a^b f(x) dx - \sum\limits_{a \leq \frac{k}{\mathfrak{n}} \leq b} p_{2k,2\mathfrak{n}} \right| 
\leq & \left| \int_a^b f(x) dx - \sum\limits_{a \leq \frac{k}{\mathfrak{n}} \leq b} d_{k,\mathfrak{n}} \right| 
	+  \left| \sum\limits_{a \leq \frac{k}{\mathfrak{n}} \leq b} d_{k,\mathfrak{n}} - 
	\sum\limits_{a \leq \frac{k}{\mathfrak{n}} \leq b} p_{2k,2\mathfrak{n}} \right|\leq\eta+\kappa. \\[8pt]  
\end{align*}
Taking  $\delta={1\over 2s}$ we obtain 
$$\eta + \kappa ={2\over 3\pi \mathfrak{n}} {2-3\delta\over (\delta(1-\delta))^{3/2}}=
{4\over 3\pi \mathfrak{n}} \left({1-{3\over 4s}}\right)\left({4s^2\over 2s-1}\right)^{{3\over 2}}={{1\over 2}C\over \mathfrak{n}}={C\over n},$$
what justifies the approximation \eqref{eq:prob_sasin} for $\delta \leq a < b \leq 1 - \delta$.  
To complete the analysis we need to investigate the errors ``on the boundaries'' of a unit interval,
\ie for $(0,\delta)$ (and, by symmetry, for $(1-\delta, \delta)$). We get
\begin{align*}
	 \left| \int_0^\delta f(x) dx - \sum\limits_{0 \leq \frac{k}{\mathfrak{n}} < \delta} p_{2k,2\mathfrak{n}} \right|
	 &  = \left| \int_0^{\frac{1}{2}} f(x) dx -  \int_\delta^{\frac{1}{2}} f(x) dx  
	 - \sum\limits_{0 \leq \frac{k}{\mathfrak{n}} \leq \frac{1}{2}} p_{2k,2\mathfrak{n}}
			    + \sum\limits_{\delta \leq \frac{k}{\mathfrak{n}} \leq \frac{1}{2}} p_{2k,2\mathfrak{n}} \right| 
			    \\
	& = \left| \frac{1}{2} -  \int_\delta^{\frac{1}{2}} f(x) dx 
	  - \frac{1}{2} + \sum\limits_{\delta \leq \frac{k}{\mathfrak{n}} \leq \frac{1}{2}} p_{2k,2\mathfrak{n}} \right|= 
	  \left|\int_\delta^{\frac{1}{2}} f(x) dx 
	  - \sum\limits_{\delta \leq \frac{k}{\mathfrak{n}} \leq \frac{1}{2}} p_{2k,2\mathfrak{n}} \right| \leq {C\over n}, \\
\end{align*}
where the last inequality follows directly from the preceding calculations. 
\qed

\medskip\par

\noindent\textbf{Remark}. 
Let $X_1, X_2, \ldots$ be  zero-average i.i.d. random variables with $E|X_i|^3<\infty$.
Denote $EX_i^2=\sigma^2$.
The central limit theorem states that $\mathcal{N}$,
a normal random variable $N(0,1)$ (denote its cdf by $\Phi$),
is the limiting distribution of 
$Y_n=\sum_{i=1}^n {X_i\over \sigma \sqrt{ n}}$ (denote its cdf by by $F^Y_n$).
It means that \textsl{for large} $n$ we can approximate $Y_n$ 
by $\mathcal{N}$ and the approximation error is bounded by the Berry-Essen inequality
$$\sup_x | F^Y_n(x) - \Phi(x)| \leq {C_0 E|X_1|^3\over \sigma^3\sqrt{n}},$$
where $C_0$ is a positive constant (in original paper  \cite{Esseen42}
it was shown that $C_0\leq 7.59$,  in \cite{Tyurin2010} 
it was shown that  $C_0\leq 0.4785$).
Lemma \ref{lem:main_error} is thus a Berry-Essen type inequality 
for approximating  $\Sasin{{n}}$ by a random variable with cdf $F^{asin}$,
tailored to our needs.

\subsubsection{Reliability of the results from the second level test}\label{sec:err_analysis_reliability}

Following \cite{Pareschi2007}, we say that a basic test (calculating $S^{asin}_{n,j}$)
is
\textsl{ not reliable} if, due to approximation errors in the computations 
of $p_j$-values, the distribution of $p_j, j=1,\ldots,m$ for 
truly random numbers is not uniform. We test the uniformity via $T^{asin}$ 
and $d^{asin}_{tv}$. Since we compare two continuous distributions, 
some discretization needs to be applied. In our testing procedure we use a partition $\mathcal{P}_s$ for this,
splitting the  interval $[0,1]$ into $(s+1)$ intervals (\ie the bins).  
Lemma \ref{lem:main_error} states that a maximum error in the 
computation of $p_j$ is bounded by ${C\over n}$ (note that $C$ implicitly depends on $s$). It means 
that a $p_j$-value that should belong to a given bin can be found in the neighboring ones
only if the distance between $p_j$ and one of the endpoints of a given bin 
is less than ${C\over n}$. Thus, this is also the fraction of $p_j$-values that 
can be found in wrong bins. The maximum propagated deviation is twice 
the error (since most bins have two neighbors), \ie 
$$\Delta = {2C\over n}.$$
Under $\mathcal{H}_0$  the distribution
of the numbers  $p_j, j=1,\ldots,m$ in the bins $1,\ldots, s+1$ is a multinomial distribution.
Indeed, this is equivalent to throwing $m$ balls independently into $s+1$ bins, where the probability of choosing first and last bin is ${1\over 2s}$
and ${1\over s}$ for all remaining bins.
The variance of the ratio of number of balls in bin $j\in\{1,s+1\}$ is equal to 
${2s-1\over m 4s^2}$,
and for bin $j\in\{2,\ldots,s\}$ is equal to ${s-1\over m s^2}$.
We have $\sigma> {\sqrt{s-1\over s^2 m}}$, where $\sigma$ is the expected 
statistical deviation of the ratio of $p_j$-values found in a given bin.
We expect that the error in approximating $p_j$-values propagates into 
an additional deviation. If the deviation is smaller than the statistical deviation,
\ie if 
\begin{equation}\label{eq:Delta_bound}
\Delta\leq \sigma,
\end{equation}
then we say that the second level test is \textsl{reliable}.
Note that the reliability of a test imposes a restriction 
on a relation between the length of a sequence used for each base test (\ie $n$) and 
the number of basic tests ($m$). Inequality (\ref{eq:Delta_bound}) implies a lower bound on $m$, namely
\begin{equation}\label{eq:m_bound}
m\leq   (s-1) \left( {n\over 2Cs}\right)^2.
\end{equation}

%
%

%
%
%

\section{The Flawed PRNG}\label{sec:flawed}
\label{subsec:experiments_flawed}



In this section we present $\textsf{Flawed}_{rng, N, \tau}$ -- a family of PRNGs.
The family depends on three parameters: $rng$ (a PRNG,
\eg the Mersenne Twister), $N$ (a small integer parameter, \eg $N = 30$) and $\tau \in [0, 1]$. $\textsf{Flawed}_{rng, N, \tau}(seed)$  generates the same output as $rng(seed)$ for
a fraction $1-\tau$ of all possible seeds. For the remaining fraction 
$\tau$ of seeds it outputs bits such that the
corresponding walk of the length $n=2^N$ spends exactly half  of the time ($2^{N-1}$ steps) above zero and exactly half of the time below zero.
In the following we will denote $n = 2\mathsf{n}=2^N$.

%

\subsection{Dyck Paths}
To generate walks with the aforementioned property we will use Dyck paths, \ie walks starting and ending at 0 with the property 
that for each prefix the number of ones is not smaller than the number of zeros.



\begin{deff}
 A sequence of $  2\mathsf{n}$  bits $B_1,\ldots, B_{2\mathsf{n}}$ is called a \textsl{Dyck path} if
 the corresponding walk $S_k$ fulfills $S_k=\sum_{i=1}^k(2B_i-1)\geq 0, 
k=1,\ldots,2\mathsf{n}-1$ and $S_{2\mathsf{n}}=0$.
 A set of all Dyck paths of length $2\mathsf{n}$ is denoted by $\mathcal{D}_{2\mathsf{n}}$.
\end{deff}
Thus, a Dyck path of length $2\mathsf{n}$ corresponds to a valid grouping of $\mathsf{n}$ pairs of 
parentheses.
We have $|\mathcal{D}_{2\mathsf{n}}|=C_\mathsf{n}={1\over \mathsf{n}+1}{2\mathsf{n}\choose \mathsf{n}}$ ($C_\mathsf{n}$ is the $\mathsf{n}$-th Catalan number).

\subsection{Sampling Dyck Paths}
We are interested in generating Dyck paths uniformly at random. To achieve this goal
we will use the following three ingredients.

\paragraph{(1) Walk sampling}
Let $\mathcal{I}_{2\mathsf{n}+1}^{-1}$ be the set of 
sequences of bits $B_1,\ldots, B_{2\mathsf{n}+1}$ such that the corresponding walk
$S_k$ ends at $-1$, \ie $S_{2\mathsf{n}+1}=-1$. 
One can easily sample a sequence $I\in\mathcal{I}^{-1}_{2\mathsf{n}+1}$ uniformly at random -- it is enough to
make a random permutation of the vector of bits $(0,\ldots,0,1,\ldots,1)$,
consisting of $\mathsf{n}+1$ zeros and $\mathsf{n}$ ones. 

\paragraph{(2) $f_{\textsf{Dyck}}$ transformation}
One can obtain a Dyck path of length $2\mathsf{n}$ from 
$I\in\mathcal{I}_{2\mathsf{n}+1}^{-1}$ using  
Algorithm~\ref{alg:fdyck}.
\begin{algorithm}[H]
\caption{$f_{\textsf{Dyck}}(I)$}\label{alg:fdyck}
\hspace*{\algorithmicindent} \textbf{Input}: $I=(B_1,\ldots,B_{2\mathsf{n}+1})\in \mathcal{I}_{2\mathsf{n}+1}^{-1}$ \\
\hspace*{\algorithmicindent} \textbf{Output}: corresponding Dyck path $f_{\textsf{Dyck}}(I)$
\begin{algorithmic}[1]
\State $S_k = \sum_{i=1}^k(2B_i-1)$ for $k = 1, \ldots, 2\mathsf{n}+1$
\State $t=\min(\textrm{argmin}\{S_k: k\in\{1,\ldots,2\mathsf{n}+1\}\})$

\State $f_{\textsf{Dyck}}(I) = (B_{t+1}, \ldots, B_{2\mathsf{n}+1}, B_1, \ldots, B_{t-1})$\\
\Return $f_{\textsf{Dyck}}(I)$
\end{algorithmic}
\end{algorithm}
\noindent
Observe that $f_{\textsf{Dyck}}$ transforms $I\in\mathcal{I}_{2\mathsf{n}+1}^{-1}$ into a Dyck path. This follows from  
simple observations:
\begin{enumerate}
 \item $I$ has exactly $\mathsf{n}+1$ zeros and $\mathsf{n}$ ones;
 \item since  $t=\min(\textrm{argmin}\{S_k: k\in\{1,\ldots,2\mathsf{n}+1\}\})$
 then $B_t = 0$ and after $B_t$ is removed then $f_{\textsf{Dyck}}(I)$ has exactly $\mathsf{n}$ bits equal to $0$ and $\mathsf{n}$ bits equal to $1$; 
 \item from the definition of $t$  (which enforces in particular that $B_{t+1} = 1$), the walk that corresponds to bits $(B_{t+1}, \ldots, B_{2\mathsf{n} + 1}, B_1, \ldots, B_{t-1})$ cannot go below $0$.
 \end{enumerate}

\begin{figure}
\centering
\begin{tikzpicture}[scale=0.9]
\begin{axis}[
    xmin=0,xmax=12.5,
    ymin=-2,ymax=3,
    grid=both,
    grid style={line width=.1pt, draw=gray!25},
    major grid style={line width=.2pt,draw=gray!25},
    axis lines=middle,
    minor tick num=1,
    enlargelimits={abs=0.5},
    axis line style={latex-latex},
    ticklabel style={font=\tiny,fill=white},
    xticklabels={a,b,2,4,6,8,10,12},
    yticklabels={a,-2,a,2,a,a},
    xlabel style={at={(ticklabel* cs:1)},anchor=north west},
    ylabel style={at={(ticklabel* cs:1)},anchor=south west}
]

 \draw[-] (axis cs: 0,0) -- (axis cs: 1,1)  ;
 \draw[-] (axis cs: 1,1) -- (axis cs: 2,0)  ;
 \draw[-] (axis cs: 2,0) -- (axis cs: 3,1)  ;
 \draw[-] (axis cs: 3,1) -- (axis cs: 4,0)  ;
 \draw[-] (axis cs: 4,0) -- (axis cs: 5,-1)  ;
  \draw[-] (axis cs: 5,-1) -- (axis cs: 6,-2)  ;
  \draw[-] (axis cs: 6,-2) -- (axis cs: 7,-1)  ;
  \draw[-] (axis cs: 7,-1) -- (axis cs: 8,-2)  ;
  \draw[-] (axis cs: 8,-2) -- (axis cs: 9,-1)  ;
  \draw[-] (axis cs: 9,-1) -- (axis cs: 10,0)  ;
  \draw[-] (axis cs: 10,0) -- (axis cs: 11,1)  ;
  \draw[-] (axis cs: 11,1) -- (axis cs: 12,0)  ;
  \draw[-] (axis cs: 12,0) -- (axis cs: 13,-1)  ;

  \draw[->,dashed, color=red!60] (axis cs: 6,-2) -- (axis cs: 6,1.5)  ;
  \draw[->,dashed, color=red!60] (axis cs: 6,-2) -- (axis cs: 12,-2)  ;

  \node at (axis cs: 6,1.7) {\footnotesize{$s$}};


%


\end{axis}
\begin{axis}[xshift=8.5cm,
    xmin=0,xmax=11.5,
    ymin=-2,ymax=3,
    grid=both,
    grid style={line width=.1pt, draw=gray!25},
    major grid style={line width=.2pt,draw=gray!25},
    axis lines=middle,
    minor tick num=1,
    enlargelimits={abs=1.0},
    axis line style={latex-latex},
    ticklabel style={font=\tiny,fill=white},
       xticklabels={a,a,2,4,6,8,10,12},
    yticklabels={a,-2,a,2,4},
    xlabel style={at={(ticklabel* cs:1)},anchor=north west},
    ylabel style={at={(ticklabel* cs:1)},anchor=south west}
]



 \draw[-] (axis cs: 0,0) -- (axis cs: 1,1)  ;
 \draw[-] (axis cs: 1,1) -- (axis cs: 2,0)  ;
 \draw[-] (axis cs: 2,0) -- (axis cs: 3,1)  ;
 \draw[-] (axis cs: 3,1) -- (axis cs: 4,2)  ;
 \draw[-] (axis cs: 4,2) -- (axis cs: 5,3)  ;
  \draw[-] (axis cs: 5,3) -- (axis cs: 6,2)  ;
    \draw[-] (axis cs: 6,2) -- (axis cs: 7,1)  ;

  \draw[-] (axis cs: 7,1) -- (axis cs: 8,2)  ;
  \draw[-] (axis cs: 8,2) -- (axis cs: 9,1)  ;
 \draw[-] (axis cs: 9,1) -- (axis cs: 10,2)  ;
 \draw[-] (axis cs: 10,2) -- (axis cs: 11,1)  ;
 \draw[-] (axis cs: 11,1) -- (axis cs: 12,0)  ;

%
%

\end{axis}

\end{tikzpicture}
\caption{A path $I\in\mathcal{I}_{13}^{-1}$ (left) and the corresponding Dyck path 
$f_{\textsf{Dyck}}(I)$ of length 
12 (right).} \label{fig:p_dyck}
\end{figure}
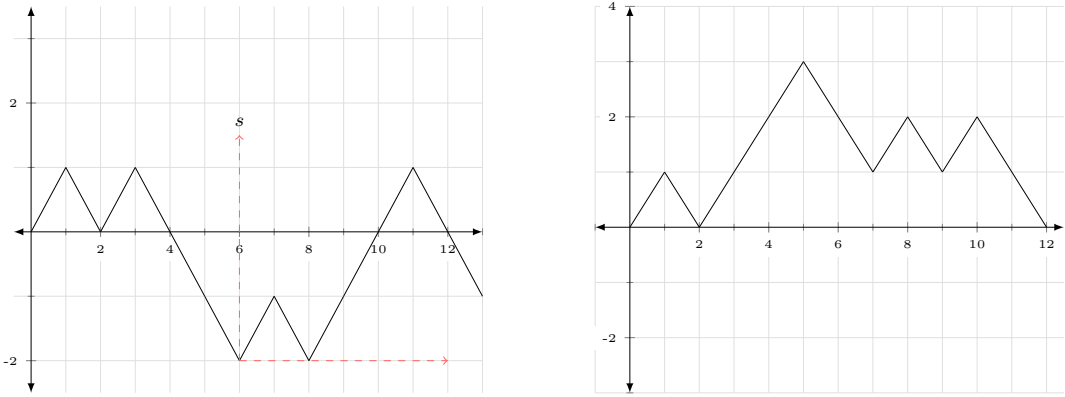

\noindent
An example of a $f_{\textsf{Dyck}}$ transformation is presented in the Figure~\ref{fig:p_dyck}.

\paragraph{(3) The Cycle Lemma} The correspondence between the set $\mathcal{I}_{2\mathsf{n}+1}^{-1}$ and the set $\mathcal{D}_{2\mathsf{n}}$ is expressed 
by the  {Cycle Lemma}~(see, \eg \cite{Dvoretzky47}).


\begin{lem}[The Cycle Lemma]\label{lem:cycle}
 For any $I\in \mathcal{I}_{2\mathsf{n}+1}^{-1}$ the path $f_{\textsf{Dyck}}(I)$ is a Dyck path.
 Moreover, any Dyck path in $\mathcal{D}_{2\mathsf{n}}$ is the image of exactly $2\mathsf{n}+1$ 
paths in
 $\mathcal{I}_{2\mathsf{n}+1}^{-1}$.
\end{lem}

\noindent 
Thus,  to obtain a  \textsl{random sampling of a Dyck path} of length $2\mathsf{n}$ one needs to run   Algorithm~\ref{alg:sample_dyck}

\begin{algorithm}[H]
\caption{$\textsf{sampleDyckPath}(\mathsf{n}, rng, seed, b)$}\label{alg:sample_dyck}
\hspace*{\algorithmicindent} \textbf{Input}:\ \  $\mathsf{n}$ -- an integer \\ 
\hspace*{\algorithmicindent} \qquad \qquad $rng$ -- a pseudorandom generator\\
\hspace*{\algorithmicindent} \qquad \qquad $seed$ -- a seed\\
\hspace*{\algorithmicindent} \qquad \qquad $b$ -- a bit deciding if the output  path should be over ($b = 0$) or under ($b = 1$) of $x$-axis\\
\hspace*{\algorithmicindent} \textbf{Output}: $P$ -- a sampled Dyck path of  length $2 \mathsf{n}$ (or its $x$-axis reflected version, for $b = 1$)

\begin{algorithmic}[1]
 \State $x = (B_1, \ldots, B_{2\mathsf{n}+1})$  where $B_1 = \ldots = B_{\mathsf{n}+1}=0$ and $B_{\mathsf{n}+2} = \ldots = B_{2\mathsf{n}+1}=1$
 \State $\sigma \leftarrow  rng(seed)$ -- a random permutation $\sigma$ of $2\mathsf{n}+1$ elements 
 \State $I = (B_{\sigma(1)}, \ldots, B_{\sigma(2\mathsf{n}+1)})$
 \State $P = f_{\textsf{Dyck}}(I)$
 \If{$b = 1$}
    \State $P = \overline{P}$  \Comment{each bit $x$ of $P$ is flipped ($x \rightarrow 1 - x$)}
 \EndIf \\
 \Return $P$
\end{algorithmic}
\end{algorithm}


%
%

\subsection{The \textsf{Flawed} generator}
As mentioned at the beginning of this section, 
the $\textsf{Flawed}_{rng, N, \tau}(seed)$  generator
--
described as Algorithm~\ref{alg:flawed} --
  works exactly the same as the underlying generator 
for a fraction $(1-\tau)$ of  seeds (lines 1-2
). For the remaining fraction of  $\tau$  of  seeds (lines (3-8), the ``else'' branch) the output is generated in the following way.
\begin{enumerate}
 \item The first $2^{N-2}$ bits are exactly the same as the first $2^{N-2}$ bits of $rng(seed)$ (line 4).
 \item The next $2^{N-2}$ bits $(z_{2^{N-2}+1}, \ldots, z_{2^{N-1}})$ are generated as follows: 
 \begin{enumerate} 
 \item a pseudorandom permutation $\pi$ is generated (line 5),
 \item the bit $z_{2^{N-2} + i}$ is set to be equal to  $1 - z_{\pi(i)}$  (lines 6-8).
 \end{enumerate} 
 As the result, there is the same number of zeros and ones in the first  $2^{N-1}$ bits -- these bits are denoted as $\mathbf{z}_L$ (\ie the corresponding walk is at zero at step $2^{N-1}$).
 
 \item The remaining $2^{N-1}$ bits 
 (denoted as $\mathbf{z}_R$) in the block are 
 obtained by calling 
 $\mathbf{z}_R \leftarrow \textsf{DyckPaths}(N, \mathbf{z}_L, rng, seed)$
 {(Algorithm \ref{alg:DyckPathsNew})}. As the result, the whole block $\mathbf{z}_L \mathbf{z}_R$  of $2^N$ output bits (concatenated blocks of $\mathbf{z}_L$ and $\mathbf{z}_R$) has the property that the corresponding walk spends the same number of steps above and below 0 (the description of $\textsf{DyckPaths}$ algorithm is below).
\end{enumerate}

\begin{algorithm}[H]
\caption{$\textsf{Flawed}_{rng, N, \tau}(seed)$}\label{alg:flawed}
\begin{algorithmic}[1]
\If{$seed \neq 0 \bmod \lceil 1/\tau \rceil$} \\
  \quad \Return $rng(seed)$
\Else 
  \State $(z_1,\ldots z_{2^{N-2}}) \leftarrow rng(seed)$
  \State $\pi \leftarrow \textsf{RandPerm}(2^{N-2}, rng(seed))$
  \For{$i = 1$ to $2^{N-2}$}
    \State $z_{2^{N-2} + i} := 1 - z_{\pi(i)}$
  \EndFor
  \State ($z_{2^{N-1}+1},\ldots, z_{2^N}) \leftarrow \textsf{DyckPaths}(N,(z_1,
\ldots, z_{2^{N-1}}), rng, seed)$ \\
  \quad \Return $(z_1, \ldots, z_{2^N})$
\EndIf
\end{algorithmic}
\end{algorithm}

%
%
\begin{exmp}[Flawed]\label{example:flawed1}
 Let  $N = 5$ and a $seed$ is such that lines 4-10 are executed. Let the result of line~4 be: $(z_1, \ldots, z_8) = (1, 1, 0, 0, 0, 0, 0, 1) \leftarrow rng(seed)$ and line~5 returns a permutation: \[\pi = \left(
    \begin{tabular}{c c c c c c c c}
    1 & 2 & 3 & 4 & 5 & 6 & 7 & 8\\
    7 & 3 & 6 & 1 & 4 & 8 & 2 & 5
    \end{tabular}\right).\]
 Then the bits computed in lines~6-8 are: \[(z_9, z_{10}, z_{11}, z_{12}, z_{13}, z_{14}, z_{15}, z_{16}) = (1 - z_7, 1 - z_3, 1 - z_6, 1 - z_1, 1 - z_4, 1 - z_8, 1 - z_2, 1 - z_5) = (1, 1, 1, 0, 1, 0, 0, 1).\]
 
 
 Then bits $z_{\mathbf{L}} = ( z_1, \ldots, z_8, z_9, \ldots, z_{16})$ are used as
 an input to function $\mathsf{DyckPaths}$ in line~10. The example is continued as Example~\ref{example:flawed2}.
\end{exmp}

Let us assume that   Algorithm \ref{alg:flawed} ($\textsf{Flawed}$) has
generated 
bits $\mathbf{z}_L=(z_1,\ldots, z_{2^{N-1}})$ by executing lines~$4-8$. Then the 
corresponding random walk spent $T_U(\mathbf{z}_L)$ time under 
 the $x$-axis and $T_O(\mathbf{z}_L)$ time
 over the $x$-axis ($T_U(\mathbf{z}_L) = 2^{N-1} - T_O(\mathbf{z}_L)$).
 The goal of the procedure $\textsf{DyckPaths}(N, \mathbf{z}_L, rnd, seed)$ is to generate bits  $\mathbf{z}_R=(z_{2^{N-1}+1}, \ldots, z_{2^{N}})
\leftarrow \textsf{DyckPaths}(N, (z_1,\ldots,z_{2^{N-1}}), rng, seed)$ 
in such a way that:
\begin{itemize}
 \item 
$T_U(\mathbf{z}_R) = T_O(\mathbf{z}_L)$ and 
\item 
$T_O(\mathbf{z}_R) = T_U(\mathbf{z}_L)$.
\end{itemize}
Then if one concatenates sequences $\mathbf{z}_L$ and $\mathbf{z}_R$, for the 
corresponding random walk it holds that $T_U(\mathbf{z}_L\mathbf{z}_R) = T_O(\mathbf{z}_L\mathbf{z}_R)$.

The following is an informal explanation of the procedure 
$\textsf{DyckPaths}(N, (z_1, \ldots, z_{2^{N-1}}), rng)$ (the formal description is provided by Algorithm~\ref{alg:DyckPathsNew}).
	\begin{enumerate}
	  \item The   walk  corresponding to $\mathbf{z}_L$ is given by $S_0:=0, S_k=\sum_{i=1}^k (2z_i-1), k=1,\ldots,2^{N-1}$. 
	  \item  The sequence $D_k$ is defined as: $D_k = \mathbbm{1}(S_k > 0 \vee S_{k-1} > 0)$, for $k = 1, \ldots, 2^{N-1}$.
	  \item The set of 
	   points where the walk changes its sign is defined as $R := \{i: D_{i+1} \neq D_i\}  \cup \{2^{N-1}\}, i=1,\ldots, 2^{N-1}-1$.

	  \item Elements of $R$ are sorted in increasing order (obtaining $(r_1, \ldots, r_w)$).
	  \item The sequence $\{l_i\}$ 
is defined as: $l_1 := 1$, and the next ``left-ends'' as $l_i = r_{i-1} + 1$ (for $i = 2, \ldots, w$).
	  \item The set  $\{O_i\}$ is defined as
 $O_i := \{ l_i, l_1+1,\ldots, r_i\}$, for $i = 1, \ldots, w$.

\item Bits $2^{N-1}+1, \ldots, 2^N$ ($\mathbf{z}_R$) are chosen  so that the 
whole walk   spends the same number of steps over and under $x$-axis.
    Dyck's paths are generated\footnote{The definitions of $l_i$ and $r_i$ imply 
    that $|O_i|$ is even,  $i=1,\ldots,w$}: 
%
   $DP_{i} = \textsf{sampleDyckPath}(|O_i|/2, rng, h(seed, D_{r_i}, i), D_{r_i})$, for $i = 1, \ldots, w$,
for some hash function $h$. We use here a hash function $h$ to obtain differently sampled Dyck's paths. This is achieved by re-seeding the generator $rng$ to be dependent on: 
	  \begin{itemize}
	  \item the current $seed$, 
	  \item a single bit equal to $D_{r_i}$ which corresponds to the type of the sequence one wants to get (over or under the $x$-axis),
	  \item the path number.
	  \end{itemize}
	  
\item A relative ordering of the paths (generated in the previous step) is obtained from a permutation 
$\rho \leftarrow  \textsf{RandPerm}(w, rng(seed))$.

    \item The resulting bits are obtained by concatenating permuted Dyck's paths.
	  
	   
	\end{enumerate}

\begin{exmp}[DyckPaths]\label{example:flawed2}
 Let input to $\mathsf{DyckPaths}$ be $z_{\mathbf{L}} = (z_1, \ldots, z_{16}) = (1, 1, 0, 0, 0, 0, 0, 1, 1, 1, 1, 0, 1, 0, 0, 1)$. 
 Then $(r_1, r_2, r_3, r_4) = (4, 10, 14, 16)$ and $(l_1, l_2, l_3, l_4) = (1, 5, 11, 15)$ and thus 
 $O_1 = \{1,2,3,4\}, O_2 = \{5,6,7,8,9,10\}, O_3 = \{11,12,13,14\}, O_4 = \{15, 16\}$.

%
 Let the output of $\textsf{sampleDyckPath}$ (called in lines 7-9 of \textsf{DyckPaths}) are $DP_1  = (0, 1, 0, 1)$, $DP_2 = (1, 1, 0, 1, 0, 0)$, $DP_3 = (0, 0, 1, 1)$, $DP_4 = (1, 0)$.
 
 Let \[\rho = \left(
    \begin{tabular}{c c c c }
    1 & 2 & 3 & 4 \\
    3 & 2 & 1 & 4 
    \end{tabular}\right).\]
    
Then $(z_{17}, \ldots, z_{32}) = DP_{\rho(1)} DP_{\rho(2)} DP_{\rho(3)} DP_{\rho(4)} = (0, 0, 1, 1, 1, 1, 0, 1, 0, 0, 0, 1, 0, 1, 1, 0)$.
\end{exmp}

\begin{algorithm}[H]
\caption{$\textsf{DyckPaths}(N, (z_1, \ldots, z_{2^{N-1}}), rng, seed)$}\label{alg:DyckPathsNew}
\begin{algorithmic}[1]
 \State $S_0:=0,\quad  S_k:=\sum_{i=1}^k (2z_i-1),\quad k=1,\ldots,2^{N-1}$ 
 \State $D_k:=\mathbbm{1}(S_k > 0 \vee S_{k-1} > 0), k=1,\ldots,2^{N-1}$
 \State $R := \{i: D_{i+1} \neq D_i\}  \cup \{2^{N-1}\}, i=1,\ldots, 2^{N-1}-1$ 
 \State Let $(r_1, \ldots, r_{w})$ be the sorted sequence of elements of $R$ 
 \State $l_i = 
 \left\lbrace
    \begin{array}{l l }
    1 & i = 1  \\
    r_{i-1} + 1 & i = 2, \ldots, w
    \end{array}\right.$

 \State $O_i := \{l_{i}, \ldots, r_{i}\}$ for $i = 1, \ldots, w$
\For{$i = 1, \ldots, w$}
    \State  $DP_i = \textsf{sampleDyckPath}(|O_i|/2, rng, h(seed, D_{r_i}, i), D_{r_i})$
 \EndFor
 \State $\rho \leftarrow  \textsf{RandPerm}(w, rng(seed))$ 
 \State $(z_{2^{N-1}+1}, \ldots, z_{2^N}) = DP_{\rho(i)} \ldots DP_{\rho(w)}$ \\
  \Return $(z_{2^{N-1}+1}, \ldots, z_{2^N})$
\end{algorithmic}
\end{algorithm}

\begin{figure}[H]
\centering
\includegraphics[scale=0.16]{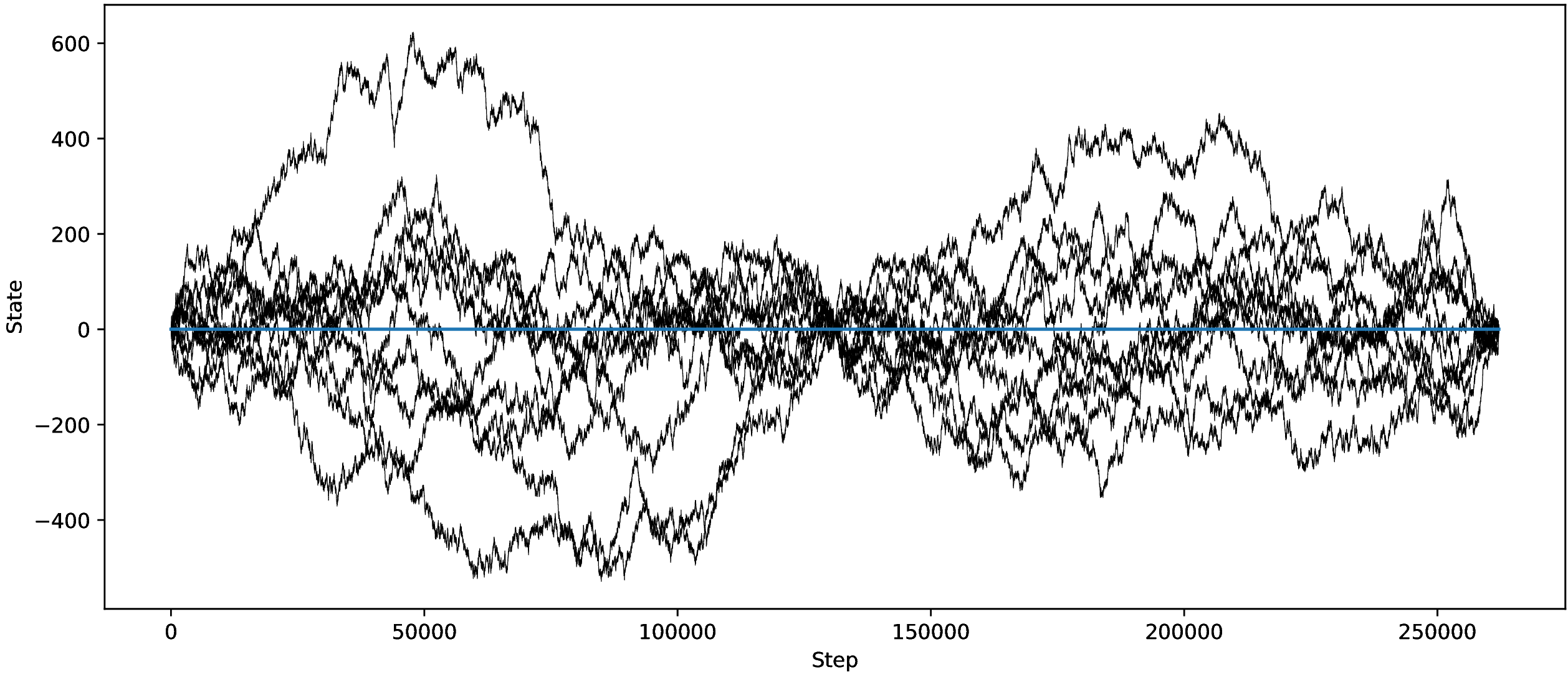}

\caption{10 trajectories of length $2^{18}$ produced by the \textsf{Flawed} generator for the seeds $\{seed_i \colon i = 0\bmod \lceil{1/\tau}\rceil\}$.}\label{fig:flaw}
\end{figure}

Ten  sample trajectories of the \textsf{Flawed} generator (all generated by the \textsl{Dyck path}-based part 
of   Algorithm \ref{alg:flawed}) are depicted in Figure~\ref{fig:flaw} (the instance of $\textsf{Flawed}_{rng, N, \tau}$ was initialized with
the following parameters: $N = {18}$, \textsl{rng} -- the Mersenne Twister). 

\section{Experimental results}\label{sec:experiments}
%
%

In this section we briefly report our experimental results of testing some
widely used PRNGs implemented in standard libraries in various programming languages.
We have applied the ASIN test to different generators including the implementations
of the standard C/C++ linear congruential generators, the standard generator \textsc{rand} from
the GNU C Library, the Mersenne Twister, the Minstd and the CMRG generator.  
As our last example we show the results of testing the $\textsf{Flawed}$ generator. $\textsf{Flawed}$ is identified
by our ASIN test as non-random, whereas it passed many other tests, including all closely related procedures
(\texttt{swalk\_RandomWalk1} test from \texttt{TestU01} with statistics: H, M, J, R, C, see Table~\ref{tab:tesu01_all}).


Each considered PRNG was tested by generating $m = 10000$ sequences of length $n \in \{2^{26},\,2^{30},\,2^{34}\}$,
using the partition $\mathcal{P}^{asin}_{40}$, \ie $s=40$. 
For these parameters our second level test is reliable (see Section \ref{sec:err_analysis_reliability}) --
$\sigma$, the expected statistical deviation of the ratio of $p_j$ values found in a given bin  is greater 
than $\sqrt{ {s-1\over s^2 m}} = 0.0015$,  what significantly exceeds the maximum propagated error  
$\Delta={2C\over n}$, 
\ie (\ref{eq:Delta_bound}) holds. Note that for $s=40$ the inequality (\ref{eq:m_bound}) yields:
\begin{itemize}
 \item $\Delta = 7.0703 \cdot 10^{-8}, m\leq 4.8760 \cdot 10^{12}$ for $n=2^{34}$, 
 \item $\Delta = 0.00000113, m\leq 1.9047 \cdot 10^{10}$ for $n = 2^{30}$, 
 \item $\Delta=0.000018$, $m\leq 7.44027\cdot 10^7$ for $n = 2^{26}$.  
\end{itemize}

%
%
%
In the experiments we used our custom implementations of tested PRNGs (except the Mersenne Twister).
We used  64-bit version of C++11 implementation of the Mersenne Twister,
\ie the class  \texttt{std::mt19937\_64}, 
which is, however, known to have some problems \cite{Harase2017}.
The generators were initialized with random seeds from \url{http://www.random.org} \cite{RandomOrg}
and each sequence was generated using different seed.

\smallskip\par

\begin{table*}[t!]
\centering
 \caption{Results of the ASIN test for  several generators with parameters $m=10000, n =2^{34}, s=40$. 
 Bolded values indicate that we reject $\mathcal{H}_0$ at the significance level $\alpha=0.0001$.}
 \label{tab:res_all}
\begin{tabular} {|l|c|c|c|c|c|c|c|c|c|c|c|}   
 \hline
   & $d^{asin}_{tv} $ & $p_{\chi^2}$ \\ \hline 
 MS Visual C++ & \textbf{0.2093} & \textbf{0.0000} \\ \hline 
 GNU C & 0.0255 & 0.2389 \\ \hline 
 Minstd 48271 & \textbf{0.2089} & \textbf{0.0000} \\ \hline 
 MT19937-64  & 0.0252 & 0.2523 \\ \hline 
%
 \end{tabular}  
\end{table*}

\smallskip\par
The results are presented in Table \ref{tab:res_all}. The values indicating that $\mathcal{H}_0$ should be rejected 
(w.r.t. significance level $\alpha=0.0001$, a value suggested by NIST for second level test)
are \textbf{bolded}. For $p_{\chi^2}$ these are simply the values smaller or equal to $\alpha$.
Concerning the values of $d^{asin}_{tv} $, Lemma \ref{lem:devroy} implies that for $\varepsilon\leq \sqrt{20\cdot 41/10000}\leq 0.2862$ we 
have $\Pro{d_{tv}^{asin}>{\varepsilon/2}}\leq 3 \exp\left(-400\varepsilon^2\right).$
It can be checked that $3exp(-400\varepsilon^2)\leq 0.0001$ for $\varepsilon\leq 0.1605$, in other words 
$$ \Pro{d_{tv}^{asin}>0.0802}\leq 0.0001,$$
\ie we reject $\mathcal{H}_0$ if the value of $d_{tv}^{asin}$ is larger than $0.0802$.
Note that for the results in Table \ref{tab:res_all} either both statistics $d_{tv}^{asin}$
and $p$-value of $T^{asin}$ reject $\mathcal{H}_0$ 
or both accept it.

%

\smallskip\par
We have also calculated the \texttt{swalk\_RandomWalk1} statistics   from \texttt{TestU01} for 
 10000 sequences of length $2^{26}$ of each  PRNG.
 The following parameters for \texttt{swalk\_RandomWalk1}  were used: $N=1, 
n=10000, r=0, s=32, L0=L1=67108864 = 2^{26}$.  The results are given in Table 
\ref{tab:tesu01_all} (including the \textsf{Flawed} generator described in Section 
\ref{subsec:experiments_flawed}).  For each Statistic H, M, J, R and C, the corresponding $p$-values were 
obtained using the
 \textsl{chi-square}
statistics.
 For convenience, $p$-values of $T^{asin}$ are  also   included in the 
Table~\ref{tab:tesu01_all} (in the column $p_{\chi^2}$).

\begin{table*}[t!]
\centering \small
 \caption{Results ($p$-values) of $T^{asin}$ and \texttt{swalk\_RandomWalk1} statistics from \texttt{TestU01} for $n=2^{26}$
(for parameters [\texttt{TestU01} notation]: $N =  1,  n = 10000,  r =  0,   s = 32,   L0 = L1 = 2^{26}$). Bolded values 
indicate that we reject $\mathcal{H}_0$ at the significance level $\alpha=0.0001$.}
 \label{tab:tesu01_all}
\begin{tabular} {|l|c|c|l|l|l|l|l|c|c|c|c|}   
 \hline
    PRNG\textbackslash Test &  $p_{\chi^2}$ 	&  Statistic  H  	
&  Statistic  M 		&  Statistic  J  	&  Statistic  R 	 
  &  Statistic  C  \\ \hline

MS Visual C++   &  0.0148 	& \quad 0.2700 & \quad 0.0900	 & 
\quad 0.6300  & \quad 0.4200 & \quad 0.8000 \\ \hline

      GNU C   &   0.4731		& \quad 0.1600 &  
\quad 0.9800	 & \quad 0.1100  & \quad 0.1900  & \quad 0.4900 \\ \hline

   Minstd 48271	  &  0.0115 			& \quad  0.0090  & 
 \quad 0.1400 & \quad 0.4900 & \quad 0.0700 &	\quad  0.0044  \\ \hline
		
   MT19937-64 	  &  0.2548 	 		& \quad 0.0800 &  
\quad 0.1000 & \quad 0.4200 & \quad 0.9700 &	\quad 0.3500 \\  \hline

$\textsf{Flawed}_{\rm MT19937-64, 26, 1/66} $ &  \textbf{0.0000}			& \quad \textbf{0.0000} &  
\quad 0.2200 & \quad \textbf{0.0000} & \quad 0.3900&	\quad 0.3800 \\  \hline

\end{tabular}  
\end{table*}


\medskip\par 
Our ASIN test would reject the MS Visual C++ PRNG and the Minstd with a multiplier 48271
(The Minstd with a multiplier 16807 gave similar results - not reported here) as 
good PRNGs. Note that this is indicated by both $p_{\chi^2}$ and the value of $d_{tv}^{asin}$. 
We also conducted the experiments 
for  the procedure \texttt{rand} from the standard library in the Borland C/C++ (not included here).
The outcomes are very akin to those for a standard PRNG in the MS Visual C++.
Note  that for  $n=2^{26}$ none of the  $p$-values calculated by the 
\texttt{swalk\_RandomWalk1} from \texttt{TestU01} suggests rejecting
the hypothesis that the MS Visual C++ PRNG is good, 
whereas the statistics H and C indicate that there can be 
some flaws in the Minstd 48271. 
It is worth mentioning that the MS Visual C++ PRNG passes the NIST Test Suite \cite{NISTtests}, as pointed out in \cite{Wang2015}.
Minstd, despite its weaknesses, became a part the C++11 standard library.
It is implemented by the classes
\texttt{std::minstd\_rand0} (with the multiplier 16807) and \texttt{std::minstd\_rand} (with the multiplier 48271).
Concerning the GNU C and the MT19937-64 -- as can be seen in both Table 
\ref{tab:res_all} and Table \ref{tab:tesu01_all} -- they can be both considered as good.
It is worth mentioning that the results  
for the CMRG generator (not reported here) were similar to those for  the MT19937-64.

The open source code of our implementation is publicly available, see 
\cite{ArcsineTest_github} (it includes the \textsf{Flawed} PRNG as well as 
  the \textsl{Law of Iterated Logarithm test} from \cite{Wang2015}).

\paragraph{Discussion on the influence of the parameter $\tau$ of the \textsf{Flawed} PRNG on the statistic $T^{asin}$}
Recall that in the Algorithm \ref{alg:flawed} the parameter $\tau$ corresponds to a fraction of simulations which
are exactly half of the time above and half of the time below $x$-axis,
 \ie we have  $S_{n,j}^{asin}=0.5$ for $\lfloor \tau m\rfloor$  simulations.
 Note that the $p_j$-value is then also equal to $0.5$.
 The remaining $\lceil(1-\tau)m\rceil$ simulations come from the $rng$. Let us assume that the $rng$ returns 
 truly random numbers. 
 
Concerning  $T^{asin}$ statistic, we have $E_1=E_{s+1}={m\over 2s}$ and $E_i={m\over s}, i=2,\ldots,s$.
 Set $r:=\lceil{s\over 2}\rceil+1$.
 For an ``ideal`` $rng$ we would have $O_1=O_{s+1}={m(1-\tau)\over 2s}, O_k={m(1-\tau)\over s}, k\in\{2,\ldots,s\}
 \setminus \{r\} $ and $O_r={m(1-\tau)\over s} +\tau m$. Thus, 
 
 $$\begin{array}{lllll}
    T^{asin}&=&\displaystyle \sum_{i=1}^{s+1} {(O_i-E_i)^2\over E_i}=(s-2) { ({ m(1-\tau)\over s} - {m\over s})^2 \over {m\over s}}
 +2 {({m(1-\tau)\over 2s} - {m\over 2s})^2\over {m\over 2s}}+ { ({m(1-\tau)\over s} +\tau m-{m\over s})^2\over 
 {m\over s}}\\[18pt]
  &=&\displaystyle (s-2){m\over s}\tau^2 + {m\over s}\tau^2 + {m\over s} (\tau(s-1))^2 = m\tau^2(s-1).\\
   \end{array}
$$
For the parameters $m=10000, \tau={1\over 66}, s=40$ we have  $T^{asin}=89.5316$ and the corresponding $p$-value is 
less than $0.0000118$. It means that even for an $rng$ producing truly random numbers, the ASIN test should recognize
the \textsf{Flawed} generator as not good.

\subsection{Results of \texttt{TestU01} for \textsf{Flawed}}

We have run several general-purpose tests against the \textsf{Flawed} generator.
For \texttt{SmallCrush} all 15 out of 15 tests were passed. For the  Mersenne Twister (MT) and the $\textsf{Flawed}_{MT, 30, 1/66}$ we run \texttt{BigCrush}. Tests for which generators failed are
presented in the Table~\ref{tabl:bigCrush}.

\vspace*{.5cm}
\begin{table}[h!]
\begin{minipage}{.5\textwidth}
\hspace*{-.5cm}\vspace*{-.5cm}
\begin{tabular}{r |  l | c | r}
tno & test name & parameters & p-value \\ \hline
74 &  RandomWalk1 R  & $L=50, r=0$  &     $6.1e-4$ \\ \hline
80 &  LinearComp  & $r = 0$ &              $1 - {\epsilon}$ \\ \hline
81 & LinearComp  & $r = 29$ &             $1 - {\epsilon}$ \\ 
\multicolumn{4}{l}{\quad} \\
\multicolumn{4}{l}{$\textsf{Mersenne Twister}$} \\ 
\end{tabular}
\end{minipage}
\begin{minipage}{.5\textwidth}
\begin{tabular}{r |  l | c | r}
tno & test name & parameters & p-value \\ \hline
80 &  LinearComp & $r = 0$  &             $1 - {\epsilon}$ \\ \hline
81 &  LinearComp & $r = 29$ &            $1 - {\epsilon}$ \\ \hline
88 &  PeriodsInStrings & $r = 0$ &         $1.1e-4$ \\ \hline
89 &  PeriodsInStrings & $r = 20$ &       $1.3e-19$ \\ \hline
102 & Run of bits & $r = 27$ &             $7.1e-4$ \\ 
\multicolumn{4}{l}{\quad} \\
\multicolumn{4}{l}{$\textsf{Flawed}_{MT, 30, 1/66}$} \\ 
\end{tabular}
\end{minipage}
\caption{Tables present which tests of the \texttt{BigCrush} were not passed. The fist column $tno$ is the test number, for the $p-value$, $\epsilon$ is value that is $\epsilon < 1.0e-15$.}
\label{tabl:bigCrush}
\end{table}

\section{Notes on Takashima's method for testing PRNGs and the arcsine test implementation from TestU01}\label{sec:takashima}
The idea of using the arcsine law for developing statistical tests for an empirical evaluation of PRNGs was formerly proposed by Takashima in \cite{Takashima94,Takashima95,Takashima}.
In this series of articles, test statistics based on the arcsine law were applied for assessing the randomness of the output of maximum-length linearly recurring sequences
($m$-sequences in short).
The experimental results presented there clearly show that the bits produced by this family of PRNGs are biased. Besides revealing the weakness of $m$-sequences,
these outcomes have also proved that Takashima's tests are effective methods, worth applying in practice. 

The approach introduced in \cite{Takashima94,Takashima95} can be briefly described as follows. After an initialization of a PRNG, 
a sequence of $2\mathfrak{n}m$ bits is generated and divided into
$m$ subsequences of length $n=2\mathfrak{n}$. Then, each subsequence  is used for constructing a random walk. 
For each of these $m$ sample random walks, the value of a test statistic based on
the arcsine law is calculated. The investigated statistic, called in \cite{Takashima94,Takashima95} the 
\textsl{sojourn time} --  denote them by $t_{n}^j, j=0,\ldots,m-1$ --  
is the time spent by a random walk above the $x$-axis.
From   $m$ realizations of  this  statistic  
an empirical distribution of 
the sojourn time   $f_{2i}=|\{j: t^j_{n}=2i\}|$, $i=0,\ldots,\mathfrak{n}$ is then derived  and compared with its theoretical 
distribution via a chi-square test.
The whole procedure is  repeated $\lambda \geq 1$ times, yielding a set of $\chi^2$ test statistics' 
values $\{\chi^2_k\},k=0,\ldots,\lambda-1$.
%
%
%
%
%
The final step of the Takashima's testing method is to count the
number of $\chi_{k}^{2}$ values falling between $90$-th and $95$-th percentile and those bigger than $95$-th 
percentile of a respective $\chi^{2}$ distribution. 
These two counts are then the basis for deciding if $\mathcal{H}_0$ should be rejected. 
Note that for $\lambda>1$ this is a third level test, which in general is \textsl{not reliable}, as shown in 
\cite{LEcuyer92}.

The author in \cite{Takashima95} considers also a slightly modified variant of the procedure, where the chi-square test is combined with the Kolomogorov-Smirnov test.
Another method, presented in \cite{Takashima}, exploits the relations between the sojourn time and the last visit time for one-dimensional random walks.

It is worth noting that in our simulations we used binary sequences of length at least $n = 2^{26}$.
Thus, a direct application of the Takashima's methods from \cite{Takashima94,Takashima95} would require large amount of additional memory to store the values of $f_{2i}, i=0,\ldots, \mathfrak{n}$.

As the arcsine law based statistical tests were proven to be useful in detecting flaws of some PRNGs,
such procedures were implemented in the \texttt{TestU01} library (see \cite{TestU01_guide}). This tool, developed by L'Ecuyer and Simard, provides a big variety of
functions for empirical examining of PRNGs. One of the test modules, \texttt{swalk}, contains a procedure \texttt{swalk\char`_RandomWalk1}, 
which calculates a bunch of test 
statistics for a sample of $m$ random walks constructed from chosen bits of generated binary sequences. 
Among them, there is the Statistic J, which implements the test based on the arcsine law. This procedure is similar to ours. Namely,
  $m$ calculated values of the test statistic
are grouped according to some partition and their empirical distribution is compared with the theoretical one by means of the chi-square test. The main difference is that in our
testing method the partition size $s$ is a parameter chosen by the user, whereas the partition used by \texttt{swalk\char`_RandomWalk1} is calculated automatically, depending on
the tested sequence. Moreover, we provide bounds on approximation errors in the computations of $p$-values (a Berry-Essen type inequality),
assuring the \textsl{reliability} of the whole testing procedure.


 \section{Conclusions}\label{sec:concl}

%
%
%
%
%
%

 In this paper we 
analyzed a method for testing PRNGs based on the arcsine law for 
random walks.  
Our procedure is a second level statistical test. We also provided a detailed error analysis of the proposed method. The approximation errors
in the calculation of $p$-values are bounded by a Berry-Essen type inequality, what allows to control the overall error, 
assuring the \textsl{reliability} of the test. We evaluate the quality of PRNGs via the chi-square statistics 
as well as by calculating a statistical distance (the total variation distance) between the empirical distribution 
of the considered characteristic for generated pseudorandom output and its theoretical distribution for truly
random binary sequences.
 
The experimental results presented in this paper show that our testing procedure can be used for detecting weaknesses
in many common PRNGs' implementations. Likewise the \textsl{Law of Iterated Logarithm test} 
from \cite{Wang2015}, the ASIN test has also revealed some flaws and
regularities in generated sequences not necessarily being identified by other current state of the art tools like the
NIST SP800-22 Testing Suite or \texttt{TestU01}. 
Thus, these kind of testing techniques seem to be very promising, as they allow also for recognition of different
kinds of deviations from those detected by existing tools. 
Nevertheless, like other statistical tests, the ASIN test is not universal and encompasses only one from an immense range of
characteristics of random bit strings and does not capture all known flaws. Therefore, the  testing procedures
relying on properties of random walks like the ASIN  test  should be used along with other tests   for more careful 
assessment of pseudorandom generators. This issue is well depicted by the provided 
example of obviously non-random generator \textsf{Flawed} for which the  LIL test has failed to detect its weaknesses, 
but the ASIN test has turned out 
to be very sensitive for that kind of deviations. 
Hence, an important line of further research is to develop another novel tests utilizing various properties of 
random walks. Such tests, when combined together, should be capable of detecting more hidden
dependencies between the consecutive bits in the sequences generated by PRNGs. This could lead to designing more robust test suites
for evaluating the quality of random numbers generated by the new implementations of PRNGs as well as those being already in use,
especially for cryptographic purposes.

 \section*{Acknowledgements}
 We would like to thank  the anonymous reviewers whose suggestions and insightful comments helped significantly 
 improve  and clarify this manuscript. In particular 
 we thank one of the reviewers for pointing out the article \cite{PareschiRS07}   on second level tests.

\section*{References}


\bibliographystyle{elsarticle-num}
\bibliography{library3}

 
\end{document}